\documentclass[11pt]{article}
\usepackage{epic,eepic,euscript,verbatim,amsthm,amsfonts,afterpage,float,bm,stmaryrd,bm,mathrsfs}
\usepackage{latexsym,amsmath,amssymb,paralist}
\usepackage{multirow}
\usepackage{comment}
\usepackage[compress]{cite}
\usepackage{algorithm}
\usepackage{algpseudocode}
\usepackage[title]{appendix}
\usepackage{xcolor}
\usepackage{BOONDOX-calo}
\usepackage{subfigure}
\usepackage{caption}
\usepackage{wrapfig}
\usepackage{algpseudocode}
\usepackage{lipsum}

\numberwithin{equation}{section}

\usepackage{graphicx}
\usepackage{epstopdf}
\usepackage{bm}
    \usepackage{nicefrac}
    \usepackage{longtable}
\usepackage{hyperref}
\usepackage[normalem]{ulem}
\newcommand{\txd}{\text d}
\newcommand{\reff}[1]{{\rm (\ref{#1})}}
\newcommand{\xm}{\text m}

\newcommand{\R}{\mathbb{R}}            
 
\newcommand{\K}{\mathbb{K}} 

\newcommand{\ve}{\varepsilon}          
\newcommand{\ep}{\epsilon}          
\newcommand{\bu}{\pmb{u}}             

\newcommand{\calM}{\mathcal M}
\newcommand{\calW}{\mathcal W}

\newcommand{\calE}{\mathcal E}

\newcommand{\cJ}{\mathcal J}

\newcommand{\calD}{\mathcal D}
\newcommand{\calA}{\mathcal A}
\newcommand{\calB}{\mathcal B}
\newcommand{\calP}{\mathcal P}
\newcommand{\calN}{\mathcal N}
\newcommand{\calL}{\mathcal L}
\newcommand{\calV}{\mathcal V}

\newcommand{\caol}{\mathcal l}

\newcommand{\x}{\mbox{\boldmath$x$}}

\newcommand{\y}{\mbox{\boldmath$y$}}

\newcommand{\n}{\mbox{\boldmath$n$}}
\newcommand{\hf}{\nicefrac{1}{2}}

\newcommand{\ciptwo}[2]{\left\langle #1 , #2 \right\rangle}

\newcommand\dt {{\Delta t}}

\newcommand{\sA}{\mathscr{A}}
\newcommand{\proofc}{{\sc Proof.}~}
\newtheorem{theorem}{Theorem}[section]
\newtheorem{lemma}[theorem]{Lemma}

\newtheorem{remark}[theorem]{Remark}

\newcommand{\td}{\tilde}

{\allowdisplaybreaks

\begin{document}

\graphicspath{{Figures/}}

\title{Entropy Increasing Numerical Methods for Prediction of Reversible and Irreversible Heating in Supercapacitors}

\author{Jie Ding\thanks{
 School of Science, Jiangnan University,  Wuxi, Jiangsu, 214122, China. E-mail: jding@jiangnan.edu.cn.}
 \and	
Xiang Ji\thanks{School of Mathematical Sciences, MOE-LSC, CMA-Shanghai, and Shanghai Center for Applied Mathematics,  Shanghai Jiao Tong University, Shanghai, 200240, China. Email: xian9ji@sjtu.edu.cn.}
\and	
Shenggao Zhou\thanks{School of Mathematical Sciences, MOE-LSC, CMA-Shanghai, and Shanghai Center for Applied Mathematics,  Shanghai Jiao Tong University, Shanghai, 200240, China. 
Shanghai Zhangjiang Institute of Mathematics, Shanghai, 201203, China. E-mail: sgzhou@sjtu.edu.cn.  }
}
\maketitle

\begin{abstract}
Accurate characterization of entropy plays a pivotal role in capturing reversible and irreversible heating in supercapacitors during charging/discharging cycles. However, numerical methods that can faithfully capture entropy variation in supercapacitors are still in lack. This work develops first-order and second-order finite-volume schemes for the prediction of non-isothermal electrokinetics in supercapacitors. Semi-implicit discretization that decouples temperature from ionic concentrations and electric potential results in an efficient first-order accurate scheme. Its numerical analysis theoretically establishes the unique solvability of the nonlinear scheme with the existence of positive ionic concentrations and temperature at discrete level.  To obtain an entropy-increasing second-order scheme, a modified Crank-Nicolson approach is proposed for discretization of the logarithm of both temperature and ionic concentrations, which is employed to enforce numerical positivity. Moreover, numerical analysis rigorously demonstrates that both first-order and second-order schemes are able to unconditionally preserve ionic mass conservation and original entropy increase for a closed, thermally insulated supercapacitor. Extensive numerical simulations show that the proposed schemes have expected accuracy and robust performance in preserving the desired properties. Temperature oscillation in the charging/discharging processes is successfully predicted, unraveling a quadratic scaling law of temperature rising slope against voltage scanning rate. It is also  demonstrated that the variation of ionic entropy contribution, which is the underlying mechanism responsible for reversible heating, is faithfully captured. Our work provides a promising tool in predicting reversible and irreversible heating in supercapacitors.


\end{abstract}

\section{Introduction}
The demand for reliable and high-performance energy storage systems becomes increasingly imperative in last decades. Supercapacitors, in which electric energy is stored and released by forming electric double layers (EDLs) through reversible ion adsorption on the solid-liquid interface of porous electrodes, have attracted tremendous attention due to their exceptional performance characteristics~\cite{FullerHarb_2018Book}.  Their unique features, such as high power density, high charging and discharging rate, long cycling life,  have shown great potential for portable electronics and hybrid vehicles~\cite{Xu_nature14, yang_sci13, Schiffer_JPS06,Ramon_JCIS18,Burke00}. 

As a classical continuum mean-field model, the Poisson--Nernst--Planck (PNP) model and its generalizations have been extensively used to describe ion transport processes or electrokinetics in EDLs~\cite{BTA:PRE:04, JiLiuLiuZhou_JPS2022}. Numerous experiments have shown that ion transport processes are accompanied by temperature oscillations during charging and discharging processes in supercapacitors~\cite{Schiffer_JPS06, Pascot_TA10, Dandeville_TA11, Zhang_TA16}. However, the PNP-type models neglect the inherently coupling between non-isothermal effects and ion transport.  In our previous work~\cite{JiLiuLiuZhou_JPS2022},  a thermodynamically consistent model, named Poisson--Nernst--Planck--Fourier (PNPF),  has been proposed to predict thermal electrokinetics in supercapacitors by using an energetic variational approach. The least action principle and maximum dissipation principle from non-equilibrium thermodynamics are utilized to derive modified Nernst--Planck equations for the description of ion transport with temperature inhomogeneity. Laws of thermodynamics are employed to develop a temperature evolution equation with heat sources arising from ionic current under electric potential differences. 

To achieve high energy storage density, electrodes with complex geometry, such as nanoporous structure, are often selected to enhance the area of EDL interfaces where the electric energy is stored.  The nanoporous electrode that hinders ion motion has profound influence on ionic electrokinetics during charging/discharging processes. 
Therefore, it is crucial to understand the impact of polydispersity and spatial arrangement of pores on the charge dynamics under non-isothermal conditions. The Transmission-Line (TL) model, a classical 1D equivalent RC circuit model, has been widely applied to probe the charging dynamics in porous electrodes~\cite{Levie_1963}. Variations of the TL model are able to account for high electric potentials~\cite{BB_PhysRevE_2010}, surface conduction in EDLs~\cite{Mirzadeh_PhysRevLett_2014}, arbitrary double-layer thickness~\cite{Gupta_PhysRevLett_2020}, or a stack-electrode model~\cite{Lian_PhysRevLett_2020}. Analogously, the heterogeneity of porous electrodes are also responsible for charge dynamics in Lithium-ion batteries. Extensive numerical and analytical studies have been conducted to understand the charge dynamics in the complex porous geometry of electrodes, contributing deeper understanding and better design strategy for energy storage systems~\cite{Schmuck_SiamAM_2015,Kirk_SiamAM_2022,Moyles_SiamAM_2019,Timms_SiamAM_2021,Vynnycky_SiamAM_2019,Fang_JComputPhys_2022}. However, not much progress has been made in capturing the geometric effect of porous electrodes on charge dynamics.


The PNP model and its generalizations possess many properties of physical significance, e.g., positivity of ionic concentration and temperature, mass conservation, and entropy increase. For instance, accurate characterization of entropy plays a crucial role in capturing reversible and irreversible heating in supercapacitors during charging/discharging cycles. Considerable efforts have been devoted to the development of numerical schemes that can maintain such properties at discrete level, ranging from finite volume schemes to finite element schemes~\cite{ProhlSchmuck09, LW14, MXL16, LW17, HuHuang_NM20, QianWangZhou_JCP21, LiuWangWiseYueZhou_2021, LiuMaimaiti_2022, DingZhou_JCP2024}. Fully implicit schemes for the PNP equations that guarantee unique solvability, unconditional positivity, and free-energy dissipation are proposed in~\cite{ShenXu_NM21, LiuWangWiseYueZhou_2021}. In the work~\cite{MXL16}, an energy-stable scheme with first-order accuracy in time for the PNP--Navier--Stokes (PNPNS) system is developed using the finite element discretization. This scheme ensures the positivity of ionic concentration by a variable transformation.  First/second-order time-stepping schemes for the PNPNS system are proposed to preserve mass conservation, positivity, and energy stability~\cite{ZhaoXu_CPC2023}. For a non-isothermal model, a first-order in time numerical scheme is proposed with guarantee of energy stability as well as positivity of the charge density and temperature via variable transformations~\cite{WuLiuZitakanov_JCP2019}. For modified PNP equations including effects arising from steric interactions and the Born solvation, both first- and second-order in time numerical schemes are developed in the works~\cite{DingWangZhou_JCP2023, DingZhou_JCP2024}, which can be proved to preserve positivity of numerical solutions, and original energy dissipation. 

In this work, we propose finite-volume schemes for the proposed PNPF model and apply them to investigate electrokinetics, temperature evolution, and entropy variation in the charging/discharging processes of supercapacitors with electrodes of complex geometry. We propose a first-order semi-implicit scheme, in which the temperature is decoupled from ionic concentrations and electric potential. The scheme is efficient to solve due to the decoupling. Rigorous numerical analysis on the proposed scheme establishes the unique existence of positive ionic concentrations and temperature at discrete level. The first-order numerical scheme is also proved to preserve mass conservation and original entropy increase in a closed,  thermally insulated system. In addition, we construct a modified Crank-Nicolson type of second-order time discretization scheme, using the logarithm of both concentrations and temperature as iteration variables. Such a second-order scheme is theoretically shown to respect mass conservation and original entropy increase as well.  Numerical studies further demonstrate that the numerical scheme has expected accuracy and presents robust performance in preserving the physical properties at discrete level. Furthermore, our model successfully predicts temperature oscillation in the charging/discharging processes, indicating that our model, together with the developed numerical methods, can robustly capture reversible and irreversible heat generation in supercapacitors. It is also found that complex electrode geometry results in an intersection point in the current-voltage loop in the cyclic voltammetry tests. Numerical investigations also unravel the temperature increasing scaling laws and entropy variation in charging/discharging phases, agreeing with existing experiments and theoretical understandings.

This paper is organized as follows. In section 2, we introduce the PNPF model and its equivalent form. In section 3, we propose a first-order scheme (Scheme I) and a second-order scheme (Scheme II) for the PNPF equations and perform numerical analysis on the proposed schemes. In section 4, we present numerical simulation results. Finally, concluding remarks are given in section 5.

\section{Non-isothermal Electrokinetic Model}
\subsection{PNPF equations}
We briefly recall the energetic variational model for the description of ionic electrodiffusion, heat generation, and thermal transport in supercapacitors~\cite{LiuWuLiu_CMS18, JiLiuLiuZhou_JPS2022}. To characterize the charging/discharging processes, we denote by $T(\x,t)$, $\psi(\x,t)$, and $c^\caol(\x,t)$ the temperature distribution, electric potential, and $\caol$th ionic concentration at location $\x$ for time $t$, respectively. Let $\Omega $ be the domain for a supercapacitor under consideration.  For any arbitrary subdomain $V\subset \Omega$, the mean-field electrostatic free-energy functional $F(V, t)$ is a functional of the particle densities and temperature given by
\begin{align}
	\label{E:freeE}
	F(V,t) = F_{pot}(V, t) + F_{ent}(V, t).
\end{align}
Here the electrostatic potential energy reads
\begin{equation}\label{ElePot}
	\begin{aligned}
		F_{pot}(V,t) &= \sum_{n,m = 1}^M\frac{q^nq^m}{2}\iint_{V}c^n(\x,t)c^m(\x^{\prime},t)G(\x,\x^{\prime})d\x d\x^{\prime} \\   
		&\qquad + \sum_{n= 1}^M q^n\int_{V}c^n(\x)\left( \psi_X(\x,t) + \int_{\Omega \backslash V} \sum_{m = 1}^M c^m(\x^{\prime}, t)G(\x,\x^{\prime})d\x^{\prime}\right) d\x,
	\end{aligned}
\end{equation}
where $q^n= z^n e$ with $z^n$ being the ionic valence and $e$ being the elementary charge. The first term represents for the Coulombic interaction energy inside $V$, with the Green function $G$ satisfying $-\nabla\cdot\epsilon_r\epsilon_0\nabla G(\x,\x^{\prime}) = \delta(\x, \x^{\prime})$. Here $\epsilon_0$ is the permittivity in vacuum and $\epsilon_r$ is the dielectric coefficient.  The second term is the electric potential energy due to the external fields, including contributions from ions outside $V$ and an external electric potential $\psi_{X}$ arising from boundary electrodes and fixed charges $\rho^f e$.
The entropy contribution in \reff{E:freeE} is given by
\[
	\begin{aligned}
		F_{ent}(V,t)  &= \int_{V} \Psi\left( c^1(\x,t),\cdots, c^M(\x,t), T(\x,t) \right)  d\x\\
							&:=\int_{V}\left[ \sum_{\caol = 1}^M\Psi_{\caol}\left( c^{\caol}(\x,t), T(\x,t)\right)  -\Psi_T(T(\x,t))\right] d\x,
	\end{aligned}
\]
where $\Psi_{\caol}$ and $\Psi_T$ are defined as
\begin{equation}\label{ent_re}
	\begin{aligned}
		\Psi_{\caol}\left(c^{\caol}(\x, t), T\left(\x, t \right)  \right)  &:= k_BT(\x,t) c^{\caol}(\x,t)\log\left( c^{\caol}(\x,t)\right) ,\\
		\Psi_T\left( T(\x, t)\right)  &:= C_Tk_BT(\x,t)\log T(\x,t).
	\end{aligned}
\end{equation}
Here $k_B$ is the Boltzmann constant and $C_T$ is a constant related to the heat capacitance of ionic species. Then the entropy $S$ is given by
\[
	S(V,t) = -\int_{V}\frac{\partial\Psi\left( c^1,\cdots ,c^M , T \right)  }{\partial T }d\x =k_B \int_{\Omega}\left( C_T (\log T +1)-\sum_{\caol=1}^M c^\caol \log c^\caol \right) d\x.
\]
To further understand contributions to entropy, we introduce the splitting $S(V,t)=S_1+S_2$ with
\begin{equation}\label{S1+S2}
S_1 = k_BC_T \int_{\Omega} (\log T +1) d\x~ \mbox{ and }~S_2 = - k_B \int_{\Omega}\sum_{\caol=1}^M c^\caol \log c^\caol d\x,
\end{equation}
where $S_1$ and $S_2$ represent contributions from temperature and ions, respectively. Application of the Legendre transform of the free energy gives the internal energy
\begin{align}
	\label{U}
	U(V,t) = F(V,t) -\int_{V}T(\x,t)\frac{\partial \Psi\left( c^1 ,\cdots, c^M , T \right)  }{\partial T }d\x.
\end{align}
Introduce the velocity $\bu^\caol$ for the $\caol$-th ionic species. Then, ionic concentration satisfies the mass conservation law
\[
\partial_t c^{\caol} +\nabla\cdot\left(c^\caol \bu^\caol \right)=0.
\]
We then employ the Least Action Principle and Maximum Dissipation Principle to derive the conservative force and dissipative force, respectively~\cite{LiuWuLiu_CMS18, JiLiuLiuZhou_JPS2022,WangLiu_Entropy_2022}. The balance of two such forces leads to 
\begin{equation}\label{fb}
	\nu^\caol c^\caol \bu^\caol =-k_B \nabla\left(c^\caol T \right)-z^\caol c^\caol e \nabla \psi,
\end{equation}
where $\nu^{\caol}$ is the viscosity of the $\caol$th species and 
$$\psi(\x, t) := \psi_X(\x, t) + \sum\limits_{\caol = 1}^Mq^{\caol}\int_{\Omega}c^{\caol}(\x^{\prime}, t) G(\x,\x^{\prime}) d\x^{\prime}$$
is the mean electric potential satisfying the Poisson's equation
\[
-\ve_0\ve_r \Delta \psi=\sum_{\caol=1}^M z^\caol  c^\caol e+\rho^fe.
\]
Combining \reff{fb} with the mass conservation law in turn leads to modified Nernst--Planck (NP) equations with non-isothermal effects:
\begin{align}
	\label{npp}
	\partial_t c^{\caol} = \nabla\cdot\frac{1}{\nu^\caol }\left[  k_B \nabla\left(c^\caol T \right)+z^\caol e c^\caol \nabla \psi\right].
\end{align}
To close the system, we next employ both first and second law of thermodynamics to  derive a temperature evolution equation
\[
C_Tk_B\partial_t T=k \Delta T-k_B T\sum_{\caol=1}^M c^\caol  \nabla\cdot\bu^\caol+\sum_{\caol=1}^M \nu^\caol c^\caol |\bu^\caol |^2, 
\]
where $k$ is the thermal conductivity and $\nu^\caol$ is the dynamical viscosity. 
The derivation details are skipped for brevity. Interested readers are referred to our previous work~\cite{JiLiuLiuZhou_JPS2022}.


To get a dimensionless formulation, we introduce the following variable rescaling
\[
\begin{aligned}
\td\x=\frac{\x}{L},~\td t=\frac{t}{\tau},~\td T=\frac{T}{T_0},~\td\psi=\frac{\psi e}{k_B T_0},~\td C=\frac{C}{c_0},~\td k=\frac{\tau k}{k_B c_0 L^2},~\td c^\caol=\frac{c^\caol}{c_0},~\td\nu^\caol=\frac{\nu^\caol}{\nu_0},\\ 
\end{aligned}
\]
where $L$ is a macroscopic length scale, $c_0$ is a characteristic concentration, $T_0$ is a characteristic temperature, $\lambda_D=\sqrt{\ve_0\ve_r k_B T_0/e^2 c_0}$ is a microscopic length scale, $\nu_0$ is a characteristic viscosity, and $\tau=\lambda_D L \nu_0/k_B T_0$ is a characteristic timescale. With the above rescaling, dropping all the tildes, we arrive at a dimensionless Poisson--Nernst--Planck--Fourier (PNPF) system
\begin{equation}\label{R_PNPF}
\left\{
\begin{aligned}
&\partial_t c^{\caol}+\ep\nabla\cdot\left(c^\caol \bu^\caol \right)=0,~~\caol=1,2,\cdots,M,\\
&\nu^\caol c^\caol \bu^\caol =- \nabla\left(c^\caol T \right)-z^\caol c^\caol \nabla \psi,~~\caol=1,2,\cdots,M,\\
&-\ep^2\Delta\psi=\sum_{\caol=1}^M z^\caol c^\caol +\rho^f,\\
&C_T\partial_t T=k \Delta T-\ep T \sum_{\caol=1}^M c^\caol  \nabla\cdot\bu^\caol+\ep\sum_{\caol=1}^M \nu^\caol c^\caol |\bu^\caol |^2,
\end{aligned}
\right.
\end{equation}
where the nondimensionalized coefficient $\ep=\frac{\lambda_D}{L}$. Consider the dimensionless boundary conditions
\begin{equation}
\left\{
\begin{aligned}
&c^\caol\bu^\caol\cdot\n=0,~\nabla T\cdot\n=0~~&&\mbox{on}~~\partial\Omega,\\
&\ep^2\frac{\partial \psi}{\partial \n}=\psi^{\rm N}~~&&\mbox{on}~~\Gamma_{\rm N},\\
&\psi=\psi^D~~&&\mbox{on}~~\Gamma_{\rm D},\\
\end{aligned}
\right.
\end{equation}
where zero-flux and thermally insulated boundary conditions are prescribed on the boundary $\partial \Omega$,  and surface charge density $\psi^{\rm N}$ and boundary potential data $\psi^D$ are prescribed on $\Gamma_{\rm N}$ and $\Gamma_{\rm D}$, with $\partial\Omega=\Gamma_{\rm N} \cup \Gamma_{\rm D}$ and $\emptyset=\Gamma_{\rm N} \cap \Gamma_{\rm D}$.
In addition, the entropy function in dimensionless form becomes
\begin{equation}\label{Sentropy}
S= \int_{\Omega}C_T (\log T +1)-\sum_{\caol=1}^M c^\caol \log c^\caol  d\x.
\end{equation}

For the derived PNPF system~\reff{R_PNPF}, one can derive a property of physical significance on entropy for a closed, thermally insulated system. 
\begin{theorem}
{\bf (Entropy Increasing)} The total entropy of a closed, thermally insulated system described by the PNPF system~\reff{R_PNPF} is increasing over time:
\begin{equation}\label{dS/dt}
\frac{dS}{dt} =\int_{V} \left( k\bigg| \frac{\nabla T}{T}\bigg|^2 + \sum_{\caol=1}^M \frac{\nu^\caol c^\caol |\bu^\caol|^2}{T} \right) d\x\geq0.
\end{equation}
\end{theorem}

\noindent \proofc
Taking derivative of $S$ with respect to time $t$, one obtains
\begin{equation}
	\begin{aligned}
		\frac{dS}{dt} &= \int_{V}\left[\frac{C_T}{T}\frac{\partial T}{\partial t} - \sum_{\caol = 1}^M\left( \log c^{\caol} + 1\right)\frac{\partial c^{\caol}}{\partial t}  \right] d\x\\
							&= \int_{V}\left[\frac{1}{T}\left( k \Delta T-\ep\sum_{\caol=1}^M c^\caol T \nabla\cdot\bu^\caol+\ep\sum_{\caol=1}^M \nu^\caol c^\caol |\bu^\caol |^2\right)  + \ep\sum_{\caol = 1}^M\left( \log c^{\caol} + 1\right)\nabla\cdot\left(  c^{\caol}\bu^{\caol} \right)   \right] d\x\\
							&=  \int_{V} \left( k\bigg| \frac{\nabla T}{T}\bigg|^2 + \sum_{\caol=1}^M \frac{\nu^\caol c^\caol |\bu^\caol|^2}{T} \right) d\x\geq0.
	\end{aligned}
\end{equation}
Here the zero-flux boundary conditions for ionic concentrations and the thermally insulated boundary condition for temperature have been used in the derivation. 
\qed

Alternatively, we can rewrite the PNPF equation into an equivalent form that facilitates the design of structure-preserving numerical methods:
\begin{equation}\label{H_PNPF}
\left\{
\begin{aligned}
&\partial_t c^{\caol}+\ep\nabla\cdot\left(c^\caol \bu^\caol \right)=0,~~\caol=1,2,\cdots,M,\\
&\nu^\caol c^\caol \bu^\caol =- c^\caol T\nabla\left(\log c^\caol \right)-c^\caol \nabla \left(z^\caol\psi + T\right),~~\caol=1,2,\cdots,M,\\
&-\ep^2\Delta\psi=\sum_{\caol=1}^M z^\caol c^\caol +\rho^f,\\
&C_T\partial_t T=k \Delta  T+T\sum_{\caol=1}^M \left[\ep\nabla\cdot\left(c^\caol \bu^\caol \log c^\caol  \right)+(1+\log c^\caol)\partial_t c^\caol  \right]+\ep\sum_{\caol=1}^M \nu^\caol c^\caol |\bu^\caol |^2. 
\end{aligned}
\right.
\end{equation}

\section{Numerical schemes}
\subsection{Finite-volume discretization and notations}
\begin{figure}[htbp]
\centering
\includegraphics[scale=.15]{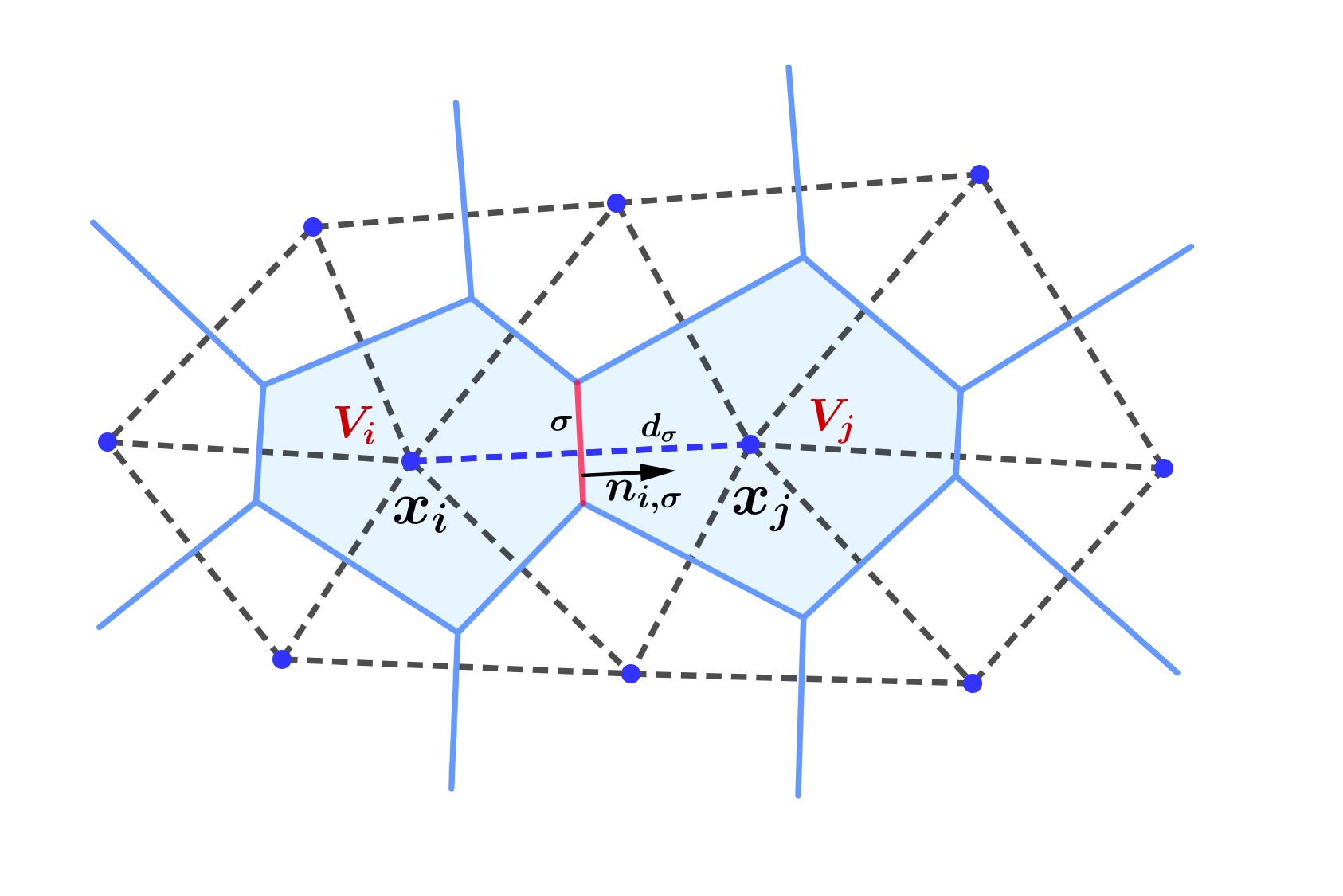}
\caption{Delaunay mesh with blue solid vertices ($\x_i$, $\x_j, \cdots$) and dual Voronoi control volumes ($V_i$, $V_j, \cdots$) with blue solid edges. The $\n_{i,\sigma}$ is the unit outward vector normal to the edge $\sigma$ of $V_i$.}
\label{f:CVs}
\end{figure}

The computational domain $\Omega\subset\R^d$ (e.g., $d=3$) is assumed to be a polygonal domain or polyhedra domain. We first introduce classical notations of finite-volume discretization mesh in the literature~\cite{Chatard2014, Chatard2012, FVM_Herbin_Book}.  The mesh $\calM=(\calV, \calE, \calP)$ covering $\Omega$ consists of a family of open polygonal control volumes $\calV:=\{V_i, i=1,2,\cdots,N \}$, a family of $d-1$ dimensional edges or faces:
 \[ 
 \begin{aligned}
 &\calE:=\calE_{int}\bigcup\calE_{ext},\\
 &\calE_{int}:=\{\sigma \subset\R^{d-1}:\sigma=\partial V_i \cap\partial V_j\},\\
 &\calE_{ext}:=\calE^{\rm D}_{ext}\cup \calE^{\rm N}_{ext}, ~\mbox{with} \left\{
 \begin{aligned}
 &\calE^{\rm D}_{ext}:=\{\sigma \subset\R^{d-1}:\sigma=\partial V_i \cap \Gamma_{\rm D}\},\\
  &\calE^{\rm N}_{ext}:=\{\sigma \subset\R^{d-1}:\sigma=\partial V_i \cap \Gamma_{\rm N}\},\\
 \end{aligned}
 \right.
 \end{aligned}
 \]
 and a family of vertices $\calP=\{\x_i, i=1,2,\cdots,N \}$; cf.~Figure~\ref{f:CVs}. The Voronoi  control volumes are defined by
 \[
V_i=\{ \y\in\Omega~\big| ~\txd(\x_i,\y)<\txd(\x_j,\y),\forall \x_j\in\calP,i\neq j\},~~i=1,2,\cdots,N,
 \]
 where $\txd(\cdot,\cdot)$ denotes the Euclidean distance in $\R^d$ and $N={\rm Card}(\calV)$. For a control volume $V_i\in\calV$, $\calE_i$ denotes the set of its edges, $\calE_{i,int}$ denotes the set of its interior edges, $\calE^{\rm D}_{i,ext}$ denotes the set of its edges included in $\Gamma_{\rm D}$, and $\calE^{\rm N}_{i,ext}$ denotes the set of its edges included in $\Gamma_{\rm N}$. 
 The size of the mesh is defined by $h=\sup\{\text{diam}(V_i),~V_i\in\calV\}$ with $\text{diam}(V_i)={\sup}_{\x,\y\in V_i} |\x-\y|$. 
Define the following three sets of indices for control volumes:
 \[
 \begin{aligned}
 &\calN_1=\{i~:~ \partial V_i\cap\calE^{\rm D}_{ext}\neq \emptyset \},~~\calN_2=\{i~:~ \partial V_i\cap\calE^{\rm N}_{ext}\neq \emptyset \},~~\calN_3=\{i~:~ \partial V_i\cap\calE_{ext}\neq \emptyset \}.\\
 \end{aligned}
 \]
For two adjacent control volumes, e.g., $V_i$ and $V_j$, the line segment $\x_i\x_j$ is orthogonal to the common edge $\sigma=\partial V_i\cap\partial V_j$. For all $\sigma\in \calE$,
 \[
 d_{\sigma}=\left\{
 \begin{aligned}
 &\txd(\x_i,\x_j)~~&&\mbox{for}~~\sigma=\partial V_i\cap\partial V_j\in \calE_{int},\\
 &\txd(\x_i,\sigma)~~&&\mbox{for}~~\sigma\in \calE_{ext}\cap\calE_i.
 \end{aligned}
 \right.
 \]
We introduce the transmissibility coefficient of an edge $\sigma\in\calE$, which is defined by $ \tau_{\sigma}=\xm(\sigma)/d_{\sigma}.$
For $\sigma\in\calE_i$, denote by $\n_{i,\sigma}$ the unit outward vector normal to $\sigma$ with respect to $V_i$. We assume that the mesh satisfies the following regularity constraint: there exists a uniform constant $C_0>0$, such that
 \[
 \txd(\x_i,\sigma)\geq C_0 \text{diam}(V_i),~~\forall V_i\in\calV,~\forall\sigma\in\calE_i.
 \]
For a function $u(\x): \Omega \to \R$, we denote by a vector of approximation values on vertices $\calP$: $u=(u_1,u_2,\cdots,u_N)^t\in\R^N,$ where the superscript $t$ represents the transpose and $u_i$ is the approximate volume average defined by
\[
u_i=\frac{1}{\xm(V_i)}\int_{V_i}u(\x)d\x ~ \mbox{ for } i =1, \cdots, N.
\]
Here $\xm(\cdot)$ denotes the measure in $\R^d$ or $\R^{d-1}$. We also label $u$ as a grid function $u: \calP \to \R $ with $u_i$ ($i=1, \cdots, N$) associated with a grid point $\x_i \in \calP$ when no confusion arises.
Define $X(\calP):= \left\{ \mbox{all grid functions } u: \calP \to \R \right\}$ and $X(\calE):= \left\{ \mbox{all edge functions } u: \calE \to \R \right\}.$
For all $u\in X(\calP)$, we define the edge function $\sA_\sigma u \in X(\calE)$ through harmonic averaging:
\[
\sA_\sigma u = \frac{[\xm(V_i)+\xm(V_j)]u_iu_j}{\xm(V_i)u_j+\xm(V_j)u_i}, \mbox{ if } \sigma=\partial V_i\cap\partial V_j \in  \calE_{int},
\]
and by corresponding boundary data with possible extrapolation if $\sigma \in  \calE_{ext}$.
For all $V_i\in\calV$, $\sigma\in\calE_i$, and $u \in X(\calP)$, we define
\[
D u_{i,\sigma}=\left\{
\begin{aligned}
&u_j-u_i~~&&\mbox{if}~~\sigma=\partial V_i\cap\partial V_j\in\calE_{i,int},\\
&u^{\rm D}_{\sigma}-u_i~~&&\mbox{if}~~\sigma\in\calE^{\rm D}_{i,ext},\\
&u^{\rm N}_{\sigma} d_\sigma~~&&\mbox{if}~~\sigma\in\calE^{\rm N}_{i,ext},
\end{aligned}
\right.
\]
where $u^{\rm D}_{\sigma}$ and $u^{\rm N}_{\sigma}$ are boundary data on the $\Gamma_{\rm D}$ and $\Gamma_{\rm N}$, respectively. For $\sigma=\partial V_i\cap\partial V_j\in\calE_{int}$ and $u \in X(\calP)$, we define 
\[
D u_\sigma= u_j-u_i,~ \mbox{with } j>i.
\]
For all $u\in X(\calP)$, we define a vector grid function $\td\nabla_h u: \calP \to \R^3 $ on a control volume $V_i$ by
\begin{equation}\label{gradtilde}
(\td\nabla_h u)_i= \left((\td D_x u)_i, (\td D_y u)_i, (\td D_z u)_i\right)^t,
\end{equation}
where
\[
(\td D_\alpha u)_i=\frac{1}{\xm(V_i)}\sum_{\sigma\in\calE_i}\xm(\sigma) \sA_\sigma  u \, n_{i,\sigma}^\alpha~ ~\mbox{for} ~ \alpha= x, y, z ~\mbox{and} ~ \n_{i,\sigma} = (n_{i,\sigma}^x, n_{i,\sigma}^y, n_{i,\sigma}^z)^t. \\
\]
For $\bu=(u,v,w)^t \in [X(\calE)]^3$, we define a grid function $\nabla_h\cdot \bu :\calP \to \R $ on a control volume $V_i$ by
\[
\left(\nabla_h\cdot \bu\right)_i=\frac{1}{\xm(V_i)}\sum_{\sigma\in\partial V_i} \xm(\sigma)\left(u_\sigma n_{i,\sigma}^x+v_\sigma n_{i,\sigma}^y +w_\sigma n_{i,\sigma}^z\right). 
\]
For grid functions $f$ and  $g \in X(\calP)$, we define
\begin{equation}\label{divu}
\nabla_h\cdot (f\nabla_h g)=\frac{1}{\xm(V_i)}\sum_{\sigma\in\calE_i}\tau_\sigma \sA_\sigma f Dg_{i,\sigma}.
\end{equation}
Then the discrete Laplacian $\Delta_h$ can be  defined by letting $f=1$. 
The discrete inner product on $X(\calP)$ is defined by
\[
\ciptwo{f}{g}=\sum_{i=1}^N \xm(V_i) f_i g_i~~~\forall f, g\in X(\calP).
\]
Then we define discrete norms for $f\in X(\calP)$ by $\|f\|_2:=\sqrt{\ciptwo{f}{f}}$ and $\|f \|_\infty:=\underset{1\leq i\leq N}{\max} |f_i|.$

We next introduce the following spaces
\[
\begin{aligned}
&\mathring{\calW}=\left\{f \in X(\calP): \ciptwo{f}{1}=0 \right\},\\
&\calW_N=\left\{f \in X(\calP): Df_{i,\sigma}=0,~ \forall  i\in\calN_3, \sigma\in\calE_{i,ext}\right\},\\
&\mathring{\calW}_N=\left\{f \in X(\calP): f\in\calW_N;~\ciptwo{f}{1}=0\right\},\\
&\calW_{\psi}=\left\{f \in X(\calP): ~f_\sigma=0,~ \forall  i\in\calN_1, \sigma\in\calE^{\rm D}_{i,ext};~Df_{i,\sigma}=0,~ \forall  i\in\calN_2, \sigma\in\calE^{\rm N}_{i,ext}\right\}.
\end{aligned}
\]
For $f^1, f^2\in\calW_N$ and $g\in X(\calP)$, inner products on edges are defined by
\[
[\nabla_h f^1,\nabla_h f^2]=\sum_{\sigma\in\calE_{int}} \tau_\sigma D f^1_{i,\sigma}D f^2_{i,\sigma} ~\mbox{ and } ~[g\nabla_h f^1,\nabla_h f^2]=\sum_{\sigma\in\calE_{int}} \tau_\sigma \sA_\sigma g\, D f^1_{i,\sigma}D f^2_{i,\sigma}.
\]
Thus, for $f\in\calW_N$, we define $\|\nabla_h f \|^2_2=[\nabla_h f,\nabla_h f].$
For $g\in\mathring{\calW}$, it is standard to verify that there exists a unique solution $\mu\in\mathring{\calW}_N$ to the equation
\[
\calL_f\mu:=-\nabla_h\cdot(f\nabla_h\mu)=g,
\]
where $f\in X(\calP)$ is positive and bounded below. This defines a linear operator $\calL_f^{-1}:\mathring{\calW} \to \mathring{\calW}_N$ by $\calL_f^{-1} g = \mu$. For $f \in X(\calP)$, $g_1, g_2\in\mathring{\calW}$, and $\calL_f\mu_l=g_l$ for $l=1, 2$, we define an inner product
\[
\ciptwo{g_1}{g_2}_{\calL_f^{-1}}:=[f\nabla_h\mu_1,\nabla_h\mu_2].
\]
Then, the discrete weighted $H^{-1}$ norm of $g$ is given by $\| g\|^2_{\calL^{-1}_f}=\ciptwo{g}{g}_{\calL^{-1}_f}$. Given $g\in X(\calP)$, it is straightforward to prove that there exists a unique solution $\phi\in \calW_{\psi}$ to the equation
\[
\calA \phi:=-\ep^2\Delta_h \phi=g.
\]
We introduce a linear operator $\calA^{-1}:X(\calP)\to \calW_{\psi}$ by $\calA^{-1}g=\phi$, . For any $g_1, g_2\in X(\calP)$, and $\calA \phi_l=g_l$ for $l=1, 2$, we define an inner product as follows:
\[
\ciptwo{g_1}{g_2}_{\calA^{-1}}:=[\nabla_h\phi_1,\nabla_h\phi_2].
\]
Then, the corresponding $H^{-1}$ norm of $g$ is given by $\| g\|^2_{\calA^{-1}}=\ciptwo{g}{g}_{\calA^{-1}}.$

The following Lemmas are straightforward and their proof is thus omitted for brevity~\cite{ChenWangWangWise_JCP2019, LiuWangWiseYueZhou_2021}.
\begin{lemma}
For $f^1, f^2\in\calW_N$, $q^1, q^2\in\calW_{\psi}$, and $g\in X(\calP)$, the following summation-by-parts formulas are valid:
\[
\begin{aligned}
&\ciptwo{f^1}{\nabla_h\cdot(g \nabla_h f^2)}=-[g\nabla_h f^1, \nabla_h f^2],\\
&\ciptwo{q^1}{\Delta_h q^2}=-[\nabla_h q^1, \nabla_h q^2].\\
\end{aligned}
\]
\end{lemma}
\begin{lemma}\label{Lem1}
Let $f\in\mathring{\calW}$, $q\in X(\calP)$, and $g\in X(\calP)$. Assume $g\geq M_0$ at each vertex in $\calP$ for some constant $M_0>0$, $ \| f \|_\infty\leq M_1$ and $ \| q \|_\infty\leq M_2$, where $M_1>0$ and $M_2>0$ may depend on $h$. The following estimates can be established:
 \[
 \|\calL^{-1}_g f \|_\infty\leq C_h:=\tilde{C}_hM_0^{-1},~~\|\calA^{-1}q \|_\infty\leq M_h,
 \]
 where $\tilde{C}_h>0$ may depend on $\Omega$, $h$, and $M_1$, and $M_h>0$ may depend on $\Omega$, $h$, and $M_2$.
 \end{lemma}

\subsection{First-order temporal discretization: Scheme I}
We now consider vertex-centered finite-volume discretization of the PNPF system~\reff{H_PNPF}. With a nonuniform time step size $\Delta t^n$ and $t^n=t^{n-1}+\Delta t^{n}$, we denote approximate solutions as grid functions $u^n\in X(\calP)$ for $u=\psi,~T,~c^\caol~(\caol=1, 2, \cdots, M)$. The discretization is based on the equivalent form~\reff{H_PNPF}. Integrating the system on each control volume and applying the divergence theorem, one obtains a vertex-centered finite volume scheme with semi-implicit discretization:
\begin{equation}\label{DisPNPF}
\left\{
\begin{aligned}
&\frac{c^{\caol,n+1}_i-c^{\caol,n}_i}{\Delta t}+ \ep\nabla_h\cdot \left( c^{\caol,n+1}\bu^{\caol,n+1}\right)_i=0,~~\caol=1,2,\cdots,M,\\
&\nu^\caol c^{\caol,n+1}\bu^{\caol,n+1}=- c^{\caol,n}\nabla_h \left(\log c^{\caol,n+1}+ z^\caol\psi^{n+1}\right)-\nabla_h\left[c^{\caol,n}(T^{n}-1)\right],\\
&-\ep^2\Delta_h\psi^{n+1}_i=\sum_{\caol=1}^M z^{\caol} c^{\caol,n+1}_i+\rho^f_i, \\
&C_T\frac{T^{n+1}_i-T^{n}_i}{\Delta t}=k\Delta_h T^{n+1}_i+ T^{n+1}_i P^{n+1}_i+\ep\sum_{\caol=1}^M \nu^\caol c^{\caol,n+1}_i |\hat{\bu}^{\caol,n+1}_i|^2,
\end{aligned}
\right.
\end{equation}
where 
\begin{equation}\label{P}
P^{n+1}_i= \sum_{\caol=1}^M \left[\ep\nabla_h\cdot ( c^{\caol,n+1}\bu^{\caol,n+1}\log c^{\caol,n+1} )_i+(1+\log c^{\caol,n+1}_i) \frac{c^{\caol,n+1}_i-c^{\caol,n}_i}{\Delta t}\right]
\end{equation}
and $\hat{\bu}^{\caol,n+1}_i$ on the control volume $V_i$ is defined by
\[
\hat{\bu}^{\caol,n+1}_i=-\frac{1}{\nu^\caol}\left[T^{n}_i\td\nabla_h\log c^{\caol,n+1}_i+\td\nabla_h\left(z^\caol\psi^{n+1}+T^n \right)_i \right].
\]
Notice that the discrete gradient in $\hat{\bu}^{\caol,n+1}_i$ follows the definition~\reff{gradtilde}, which is defined via integration by parts in the control volume $V_i$. We refer to the scheme \reff{DisPNPF} as ``{\bf Scheme I}" in what follows. Notice that in Scheme I, the ionic concentration and electric potential are coupled and can be solved as in the classical PNP schemes~\cite{LiuWangWiseYueZhou_2021}. With obtained $c^{\caol,n+1}$ and $\psi^{n+1}$ in $P^{n+1}_i$, the temperature equation is decoupled and can be efficiently solved with a linear scheme.  

The initial and boundary conditions are discretized as follows:
\begin{equation}\label{Bcs}
\begin{aligned}
&c^{\caol,0}_i=\frac{1}{\xm(V_i)}\int_{V_i} c^{\caol,0}(\x)d\x,~~&&\caol=1,2,\cdots,M,~~\forall V_i\in\calV,\\
&c^{\caol,n+1}_{\sigma}\bu^{\caol,n+1}_{\sigma}\cdot\n=0,~D T^{n+1}_{i,\sigma}=0,~~~&&\forall ~i\in\calN_3,~~\sigma\in\calE_{i,ext}=\calE^{\rm D}_{i,ext}\cup \calE^{\rm N}_{i,ext},\\
&\psi^{\rm D}_{\sigma}=\frac{1}{\xm(\sigma)}\int_{\sigma} \psi^{\rm D}(l)\,d l~~&& \forall~ i\in\calN_1,~~\sigma\in\calE^{\rm D}_{i,ext},
\\
&\psi^{\rm N}_{\sigma}=\ep^2 D\psi^n_{i,\sigma}/d_{\sigma},~~&&\forall ~i\in\calN_2,~~\sigma\in\calE^{\rm N}_{i,ext}.
\end{aligned}
\end{equation}
We next prove several properties of physical significance that are preserved by Scheme I at discrete level in the following theorems.

\begin{theorem}\label{th:conservation}
{\bf (Mass conservation)} The Scheme I \reff{DisPNPF} respects the mass conservation law:
\[
\ciptwo{c^{\caol,n+1}}{1}=\ciptwo{c^{\caol,n}}{1}~~\mbox{for }~\caol=1,2,\cdots,M.
\]
\end{theorem}
\noindent \proofc
It follows from \reff{divu} and the Scheme I \reff{DisPNPF} that
\[
\begin{aligned}
\sum_{i=1}^N \xm(V_i) \left(c^{\caol,n+1}_{i}-c^{\caol,n}_{i}\right)=-\ep\Delta t\sum_{i\in\calN_3}\sum_{\sigma\in\calE_{i,ext}} \xm(\sigma)c^{\caol,n+1}_{\sigma}\bu^{\caol,n+1}_{\sigma}\cdot \n=0,~~\caol=1,2,\cdots,M,
\end{aligned}
\]
where the zero-flux boundary conditions \reff{Bcs} have been used. This completes the proof. 
\qed

 \begin{theorem}\label{th:EUP}
{\bf (Existence, uniqueness, and positivity)}
Let $C_{\rm min}^{\caol,n}={\min}_{1\leq i\leq N}c^{\caol,n}_{i}$. If $C_{\rm min}^{\caol,n}>0$ and $\|c^{\caol,n}\|_\infty <\infty$, then there exists unconditionally a unique solution of the electric potential $\psi^{n+1}$ and ionic concentrations $c^{\caol,n+1}$ ($\caol = 1 , \cdots, M$) to the Scheme I \reff{DisPNPF}, such that $c^{\caol,n+1}_{i}>0$ for $i=1,\cdots,N$. If, in addition, $T_i^{n}>0$ ($i=1,\cdots,N$) and the time step size $\dt < C_T/C_*$ with $C_*:= \| P^{n+1} \|_\infty$ dependent on $\psi^{n+1}$ and $c^{\caol,n+1}$, then there exists a unique solution of temperature $T^{n+1}$ to the Scheme I \reff{DisPNPF} such that
$
T^{n+1}_i>0 ~~\mbox{for}~~ i=1,2,\cdots,N.
$
\end{theorem}

\noindent \proofc First, we prove that there is a unique positive numerical solution of ionic concentrations and electric potential to the nonlinear scheme~\reff{DisPNPF}. We construct an auxiliary variable $\beta^{\caol,n}$ that satisfies
\[
\left\{
\begin{aligned}
&-\nabla_h\cdot \left[c^{\caol,n}\nabla_h\beta^{\caol,n} \right]_i=\Delta_h \left(c^{\caol,n}(T^n-1) \right)_i,~~ &&i=1,2,\cdots,N,\\
&c^{\caol,n}D\left( \beta^{\caol,n} \right)_{i,\sigma}=-D\left( c^{\caol,n}(T^n-1) \right)_{i,\sigma},~~ &&\forall ~i\in\calN_3,~\sigma\in\calE_{i,ext}.
\end{aligned}
\right.
\]
It can be readily checked that the source term and the Neumann boundary data of this equation satisfy the compatibility condition. Hence, it is solvable and has a unique solution $\beta^{\caol,n} \in \mathring{\calW}$. Notice that $\beta^{\caol,n}$ is independent of $c^{\caol,n+1}$ and $\psi^{n+1}$.  With the constructed $\beta^{\caol,n}$, the discrete equation for concentrations can be converted to
\begin{equation}\label{muEq}
\left\{
\begin{aligned}
&-\nabla_h\cdot c^{\caol,n}\nabla_h\left(\log c^{\caol,n+1}+z^\caol\psi^{n+1}-\beta^{\caol,n} \right)_i=-\frac{\nu^\caol (c^{\caol,n+1}-c^{\caol,n})_i}{\ep\dt}, &&i=1,2,\cdots,N,\\
&D\left(\log c^{\caol,n+1}+z^\caol\psi^{n+1}-\beta^{\caol,n}  \right)_{i,\sigma}=0, &&\forall ~i\in\calN_3,~\sigma\in\calE_{i,ext}.\\
\end{aligned}
\right.
\end{equation}
It follows from Theorem \ref{th:conservation} on ionic mass conservation that $c^{\caol,n+1}-c^{\caol,n} \in \mathring{\calW}$, which further implies that the equation \reff{muEq} satisfies the compatibility condition and has a unique solution $\log c^{\caol,n+1}+z^\caol\psi^{n+1}-\beta^{\caol,n} \in \mathring{\calW}_N$. Thus, we have $\log c^{\caol,n+1}+z^\caol\psi^{n+1}=\beta^{\caol,n}- \frac{\nu}{\ep \dt}\calL^{-1}_{ c^{\caol,n}}(c^{\caol,n+1}-c^{\caol,n}) \in \mathring{\calW}$.

In addition, we introduce a splitting of the electric potential: $\psi^{n+1}= \psi_1^{n+1}+ \psi_2^{n+1}$, where $\psi_1^{n+1}$ and  $\psi_2^{n+1}$ satisfy 
\[
\begin{aligned}
\left\{
\begin{aligned}
&-\ep^2\Delta_h\psi^{n+1}_{1,i}=\sum_{\caol=1}^M z^\caol c^{\caol,n+1}_i, &&i=1,2,\cdots,N,\\
&\psi^{n+1}_{1,\sigma}=0, &&\forall~ i\in\calN_1,~\sigma\in\calE^{\rm D}_{i,ext},\\
&\ep^2D\psi^{n+1}_{1,i,\sigma}=0, &&\forall ~i\in\calN_2,~\sigma\in\calE^{\rm N}_{i,ext},\\
\end{aligned}
\right.
\quad
\left\{
\begin{aligned}
&-\ep^2\Delta_h\psi^{n+1}_{2,i}=\rho^f_i, &&i=1,2,\cdots,N,\\
&\psi^{n+1}_{2,\sigma}=\psi^{\rm D}_\sigma, &&\forall~ i\in\calN_1,~\sigma\in\calE^{\rm D}_{i,ext},\\
&\ep^2D\psi^{n+1}_{2,i,\sigma}/d_\sigma=\psi^{\rm N}_\sigma, &&\forall ~i\in\calN_2,~\sigma\in\calE^{\rm N}_{i,ext},
\end{aligned}
\right.
\end{aligned}
\]
respectively. Notice that $\psi_2^{n+1}$ is independent of $c^{\caol,n+1}$ but could depend on time, since the boundary data, e.g., $\psi^{\rm D}$, could depend on time. By the equation for $\psi_1^{n+1}$, we have $\psi_1^{n+1}=\calA^{-1}(\sum_{\caol=1}^M z^\caol c^{\caol,n+1}) \in \calW_{\psi}.$ 

It can be verified that $c^{\caol,n+1}$ and $\psi^{n+1}$ correspond to the minimizer of a constructed discrete energy functional
\[
\begin{aligned}
\cJ(c)=&\sum_{\caol=1}^M\left[\frac{\nu^\caol}{2\ep\Delta t} \left \|c^\caol-c^{\caol,n} \right \|^2_{\calL^{-1}_{ c^{\caol,n}}}+\ciptwo{c^\caol}{\log c^\caol -1}+\ciptwo{c^\caol}{ z^\caol \psi_2^{n+1}-\beta^{\caol,n}}\right] +\hf \left\| \sum_{\caol=1}^M z^\caol c^\caol \right\|^2_{\calA^{-1}}  
\end{aligned}
\]
over the admissible set
\[
\K_{h}:=\left\{c\bigg| 0 <c^\caol_{i}< \xi_\caol,~\ciptwo{c^{\caol}}{1}=\alpha_\caol \xm(\Omega),~\mbox{for}~ \caol=1,\cdots,M,~i=1,\cdots,N \right\},
\]
where $c=(c^1,c^2,\cdots,c^M)$, $\alpha_\caol =\ciptwo{c^{\caol,0}}{1}/\xm(\Omega)$ is the mean concentration of the $\caol$th ionic species, and $\xi_\caol=\frac{\alpha_\caol}{\min_i \{\xm(V_i)\}}$.
Consider a closed set $\K_{h,\delta}\subset\K_h$:
\[
\K_{h,\delta}:=\left\{c\bigg| \delta \leq c^\caol_{i}\leq \xi_\caol-\delta,~ \ciptwo{c^{\caol}}{1}=\alpha_\caol \xm(\Omega),~\mbox{for}~ \caol=1,\cdots,M,~i=1,\cdots,N \right\},
\]
where $\delta\in(0,\frac{\xi_\caol}{2})$. Obviously, $\K_{h,\delta}$ is a bounded, convex, and compact subset of $\K_h$. By the convexity of $\cJ$, there exists a unique minimizer of $\cJ$ in $\K_{h,\delta}$. 

 Suppose the minimizer of $\cJ$ over $\K_{h,\delta}$, $c^*=(c^{1,*},c^{2,*},\cdots,c^{M,*})$, touched the boundary.
Without loss of generality, we assume that there exists a control volume $V_{i_0}$ and the $\caol_0$th ionic species such that $c^{\caol_0,*}_{i_0}=\delta$.  Assume the maximum of $c^{\caol_0,*}$ is achieved at the control volume $V_{i_1}$. 
Clearly, 
\[
c^{\caol_0,*}_{i_1}\geq \alpha_{\caol_0},~~c^{\caol_0,*}_{i_0}\leq \alpha_{\caol_0}.
\]
Let 
\[
d^{\caol_0,*}_t=c^{\caol_0,*}+t e_{i_0}-t e_{i_1},
\]
where $t\in \R^+$ and $e_{i}\in X(\calP)$ is a grid function being $1$ at the grid point $\x_i$ and $0$ at the rest. Denote $d^{*}_t= (c^{1,*}, \cdots, d^{\caol_0,*}_t, \cdots, c^{M,*})$. Consider sufficiently small $t>0$, then $d^{*}_t\in \K_{h,\delta}$. 
Direct computation yields
\begin{equation}\label{J:eq1}
\begin{aligned}
\underset{t\rightarrow 0+}{\lim}\frac{\cJ(d_t^{*})-\cJ(c^{*})}{t}
=&\frac{\nu^{\caol_0}}{\ep\Delta t}\left[\xm(V_{i_0})\calL^{-1}_{ c^{\caol_0,n}}\left(c^{\caol_0,*}-c^{\caol_0,n} \right)_{i_0}-\xm(V_{i_1})\calL^{-1}_{ c^{\caol_0,n}}\left(c^{\caol_0,*}-c^{\caol_0,n} \right)_{i_1}\right]\\
&+z^{\caol_0} \xm(V_{i_0})\calA^{-1} \left(\sum_{\caol=1}^M z^\caol c^\caol\right)_{i_0}-z^{\caol_0}\xm(V_{i_1})\calA^{-1} \left(\sum_{\caol=1}^M z^\caol c^\caol\right)_{i_1}\\
&+\xm(V_{i_0})\left(z^{\caol_0}\psi^{n+1}_{2,i_0} -\beta^{\caol,n}_{i_0}\right)-\xm(V_{i_1})\left(z^{\caol_0}\psi^{n+1}_{2,i_1}-\beta^{\caol,n}_{i_1}\right)\\
&+\xm(V_{i_0})\log c^{\caol_0,*}_{i_0} -\xm(V_{i_1})\log c^{\caol_0,*}_{i_1}.
\end{aligned}
\end{equation}
Since $\|c^{\caol_0,*}-c^{\caol_0,n}  \|_\infty\leq \xi_{\caol_0}+\|c^{\caol,n}\|_\infty$, $c^{\caol_0,n}\geq C_{\rm  min}^{\caol_0,n}>0$, and $\|\sum_{\caol=1}^M z^\caol c^\caol \|_\infty\leq \sum_{\caol=1}^M |z^\caol|\xi_\caol$, we have by Lemma \ref{Lem1} that 
\begin{equation}\label{J:eq2}
\begin{aligned}
&\bigg\|\xm(V_{i_0})\calL^{-1}_{ c^{\caol_0,n}}\left(c^{\caol_0,*}-c^{\caol_0,n} \right)_{i_0}-\xm(V_{i_1})\calL^{-1}_{ c^{\caol_0,n}}\left(c^{\caol_0*}-c^{\caol_0,n} \right)_{i_1}\bigg\|_\infty\leq C_h,\\
&\bigg\| z^{\caol_0}\xm(V_{i_0})\calA^{-1} \left(\sum_{\caol=1}^M z^\caol c^\caol\right)_{i_0}-z^{\caol_0}\xm(V_{i_1})\calA^{-1} \left(\sum_{\caol=1}^M z^\caol c^\caol\right)_{i_1}\bigg\|_\infty\leq M_h,
\end{aligned}
\end{equation}
where $C_h$ and $M_h$ may depend on $C_{\rm  min}^{\caol_0,n}$, $\|c^{\caol,n}\|_\infty$, and $h$. Furthermore, it is readily seen that
\begin{equation}\label{J:eq3}
\begin{aligned}
\bigg\| \xm(V_{i_0})\left(z^{\caol_0}\psi^{n+1}_{2,i_0} -\beta^{\caol,n}_{i_0}\right)-\xm(V_{i_1})\left(z^{\caol_0}\psi^{n+1}_{2,i_1}-\beta^{\caol,n}_{i_1}\right) \bigg\|_\infty\leq C^*,
\end{aligned}
\end{equation}
where $C^*>0$ is independent of $\delta$. 
It follows from $c^{\caol_0,*}_{i_0}=\delta$ and $c^{\caol_0,*}_{i_1}\geq \alpha_{\caol_0}$ that
\begin{equation}\label{J:eq4}
\xm(V_{i_0})\log c^{\caol_0,*}_{i_0}-\xm(V_{i_1})\log c^{\caol_0,*}_{i_1}\leq \xm(V_{i_0})\log\delta-\xm(V_{i_1})\log\alpha_{\caol_0}.
\end{equation}
Combination of \reff{J:eq1}-\reff{J:eq4} yields
\[
\begin{aligned}
\underset{t\rightarrow 0+}{\lim}\frac{\cJ(d_t^{*})-\cJ(c^{*})}{t}
\leq&\xm(V_{i_0})\log\delta-\xm(V_{i_1})\log\alpha_{\caol_0}+\frac{\nu^{\caol_0}}{\ep\Delta t}C_h+M_h+C^*.
\end{aligned}
\]
Clearly, as $\delta \rightarrow 0+ $, 
\[
\xm(V_{i_0})\log\delta-\xm(V_{i_1})\log\alpha_{\caol_0}+\frac{\nu^{\caol_0}}{\ep\Delta t}C_h+M_h+C^*<0.
\]
Thus, we arrive at
\[
\cJ(d_t^{*})<\cJ(c^{*}), ~ \mbox{as~} \delta \rightarrow 0+. 
\]
This contradicts the assumption that, when $\delta$ is sufficiently small,  the minimizer of 
$\cJ$, $c^{*}\in \K_{h,\delta}$, touched the boundary of $\K_h$. 
Therefore, the minimizer of $\cJ$ is achieved at an interior point, i.e., $c^{*}\in\mathring{\K}_{h,\delta} \subset \K_h$ as $\delta \rightarrow 0$. This establishes the desired existence of $c^{\caol,n+1}_{i}>0$ for $i=1,2,\cdots,N,$ and $\caol = 1, \cdots, M.$ In addition, the strict convexity of $\cJ$ over $\K_h$ implies the uniqueness of the numerical solution. The unique existence of $c^{\caol,n+1}$ directly implies that of $\psi^{n+1}$ through the discrete Poisson equation.  

We now proceed to prove the unique existence of positive temperature. The temperature equation is linear
and can be rewritten in a matrix form
\[
\left(\frac{C_T}{\dt}\calE-k\calB-\calD \right)T^{n+1}=\frac{C_T}{\dt}T^n+\theta,
\]
where $T^{n+1}$ and $T^n$ are column vectors with entries $T^{n+1}_i$ and $T^n_i$, respectively, $\calE$ is the identity matrix, $\calB$ corresponds to the discrete Laplacian operator $\Delta_h$, $\calD$ is a diagonal matrix with entries given by $P_i^{n+1}$, 
and $\theta$ is a positive column vector with entries given by $\ep\sum_{\caol=1}^M \nu^\caol c_i^{\caol,n+1} |\hat{\bu}_i^{\caol,n+1}|^2$ for $i=1, \cdots, M$. By $\|\calD\|_{\infty}\leq C_*$ and $\dt < C_T/C_*$,
 we obtain that $\left(\frac{C_T}{\dt}\calE-k\calB-\calD \right)$ is an M-matrix with $\left(\frac{C_T}{\dt}\calE-k\calB-\calD \right)^{-1}>0$ in the element-wise sense. Thus, since $T^{n}>0$, the desired unique existence of positive temperature has been established. 
\qed


The discrete entropy is defined by
\begin{equation}\label{DisEntropy}
S^n_h=-\sum_{\caol=1}^M \ciptwo{c^{\caol,n}}{\log c^{\caol,n}}+\ciptwo{\log T^n+1 }{C_T},
\end{equation}
which is second-order spatial discretization of the continuous entropy functional~\reff{Sentropy}.  
\begin{theorem}\label{th:entropy}
{\bf (Discrete entropy increasing)}
The solution to the Scheme I \reff{DisPNPF} retains original entropy increasing, i.e.,
\begin{equation}\label{ED1:eq0}
\frac{S^{n+1}_h-S^n_h}{\dt}\geq\ep\sum_{\caol=1}^M \ciptwo{\nu^\caol c^{\caol,n+1} |\hat{\bu}^{\caol,n+1}|^2}{\frac{1}{T^{n+1}}}- k \sum_{\sigma\in\calE_{int}} \tau_\sigma D T^{n+1}_{i,\sigma} D\left(\frac{1}{T^{n+1}}\right)_{i,\sigma}\geq 0.
\end{equation}
\end{theorem}

\noindent \proofc 
Taking a discrete inner product of the fourth equation in \reff{DisPNPF} with $\frac{\dt }{T^{n+1}}$ leads to
\[
\begin{aligned}
&C_T\ciptwo{T^{n+1}-T^{n}}{ \frac{1}{T^{n+1}}}\\
&=\dt k \ciptwo{\Delta_h T^{n+1}}{\frac{1}{T^{n+1}}}+\dt\ep\sum_{\caol=1}^M \ciptwo{\nu^\caol c^{\caol,n+1} |\hat{\bu}^{\caol,n+1}|^2}{\frac{1}{T^{n+1}}}\\
&\qquad+\dt \ep\sum_{\caol=1}^M \left[\ciptwo{\nabla_h\cdot (c^{n+1}\bu^{\caol,n+1} \log c^{\caol,n+1}) }{1}+\ciptwo{(1+\log c^{\caol,n+1}) \frac{c^{\caol,n+1}-c^{\caol,n}}{\ep\Delta t}}{1}\right],
\end{aligned}
\]
where the equation \reff{P} has been used. For the first term on the right hand side, one has
\[
\begin{aligned}
\ciptwo{\Delta_h T^{n+1}}{\frac{1}{T^{n+1}}}=&\sum_{i=1}^N \frac{1}{T^{n+1}_i}\sum_{\sigma\in\partial V_i} \tau_\sigma D T^{n+1}_{i,\sigma}=-\sum_{\sigma\in\calE_{int}} \tau_\sigma D T^{n+1}_{i,\sigma} D\left(\frac{1}{T^{n+1}}\right)_{i,\sigma} \geq 0,
\end{aligned}
\]
where the thermally insulated boundary condition~\reff{Bcs} for $T$ has been used in the second equality, and monotonicity of the function $\frac{1}{x}$ has been used in the last inequality.
The following inequalities are available for the inner products:
\begin{equation}\label{ED1:eq3}
\begin{aligned}
&\ciptwo{T^{n+1}-T^{n}}{ \frac{1}{T^{n+1}}}=\ciptwo{1-\frac{T^n}{T^{n+1}}}{1}\leq \ciptwo{\log T^{n+1}-\log T^n}{1},\\
&\ciptwo{\nabla_h\cdot (c^{n+1}\bu^{\caol,n+1} \log c^{\caol,n+1}) }{1}=0,\\
&\ciptwo{\nu^\caol c^{\caol,n+1} |\hat{\bu}^{\caol,n+1}|^2}{\frac{1}{T^{n+1}}}\geq 0.\\\end{aligned}
\end{equation}
Applying the inequality $ y\log y-x\log x\geq (1+\log y)(y-x)$ with $y=c^{\caol,n+1}_i$ and $x=c^{\caol,n}_i$,  we find
\[
\ciptwo{(1+\log c^{\caol,n+1}) \frac{c^{\caol,n+1}-c^{\caol,n}}{\ep\Delta t}}{1}
\geq\frac{1}{\dt}\left( \ciptwo{c^{\caol,n+1} }{ \log c^{\caol,n+1}}-\ciptwo{c^{\caol,n}}{ \log c^{\caol,n}} \right).
\]
 The proof of \reff{ED1:eq0} is completed by combining the above results. \qed

\subsection{Second-order temporal discretization: Scheme II}
Straightforward second-order extension of the above temporal discretization may lead to one that fails to respect the entropy increasing property~\reff{dS/dt} at discrete level. Special care should be taken when dealing with $\log T$ and $\sum_{\caol=1}^M c^\caol \log c^\caol$ in the entropy functional~\reff{Sentropy}.  In order to achieve second-order temporal accuracy while still preserving the entropy increasing and positivity property, we take  $\xi = \log T$ and $\eta^\caol = \log c^\caol$ for $\caol=1, 2, \cdots, M$ as working variables, and propose a novel modified Crank-Nicolson discretization for the PNPF system~\reff{R_PNPF} as follows: 
\begin{equation}\label{DisPNPF2}
\left\{
\begin{aligned}
&\frac{e^{\eta^{\caol,n+1}_i}-e^{\eta^{\caol,n}_i}}{\Delta t}+ \ep\nabla_h\cdot (e^{\eta^{\caol,n+\hf}}\bu^{\caol,n+\hf})_i=0,~~\caol=1,2,\cdots,M,\\
&\nu^\caol e^{\eta^{\caol,n+\hf}}\bu^{\caol,n+\hf}=-e^{\eta^{\caol,n+\hf}+\xi^{n+\hf}}  \nabla_h Q^{\caol,n+\hf}
-e^{\eta^{\caol,n+\hf}}\nabla_h(z^\caol\psi^{n+\hf}+e^{\xi^{n+\hf}}) ,\\
&-\ep^2 \Delta_h \psi^{n+1}_i=\sum_{\caol=1}^M z^{\caol}e^{\eta^{\caol,n+1}_i}+\rho^f_i, \\
&C_T\frac{e^{\xi^{n+1}_i}-e^{\xi^{n}_i}}{\Delta t}=-k\nabla_h\cdot e^{\xi^{n+\hf}} \nabla_h(\log R^{n+\hf})_i+\frac{P^{n+\hf}_i}{R^{n+\hf}_i}+\ep\sum_{\caol=1}^M \nu^\caol e^{\eta^{\caol,n+\hf}_i} |\check{\bu}^{\caol,n+\hf}_i|^2,
\end{aligned}
\right.
\end{equation}
where 
$Q^{\caol,n+\hf}$, $P^{n+\hf}$, $R^{n+\hf}$, $\psi^{n+\hf}$, and $e^{f^{n+\hf}}$ are given by
\begin{equation}\label{PQR}
\begin{aligned}
Q^{\caol,n+\hf}=&\eta^{\caol,n+1}-\frac{1}{2e^{\eta^{\caol,n+1}}}(e^{\eta^{\caol,n+1}}-e^{\eta^{\caol,n}})-\frac{1}{6e^{2\eta^{\caol,n+1}}}(e^{\eta^{\caol,n+1}}-e^{\eta^{\caol,n}})^2,\\
P^{n+\hf}=&\sum_{\caol=1}^M \left[\ep\nabla_h\cdot (Q^{\caol,n+\hf}e^{\eta^{\caol,n+\hf}}\bu^{\caol,n+\hf} )+(1+Q^{\caol,n+\hf}) \frac{e^{\eta^{\caol,n+1}}-e^{\eta^{\caol,n}}}{\Delta t}\right],\\
R^{n+\hf}=&e^{-\xi^{n+1}}+\frac{e^{\xi^{n+1}}-e^{\xi^n}}{2e^{2\xi^{n+1}}}+\frac{(e^{\xi^{n+1}}-e^{\xi^n})^2}{3e^{3\xi^{n+1}}},\\
\psi^{n+\hf}=&\frac{1}{2}(\psi^n+\psi^{n+1}),\\
e^{f^{n+\hf}}=&e^{\frac{1}{2}(f^{n}+f^{n+1})} \mbox{ for } f=\xi,~\eta^\caol,\\
\end{aligned}
\end{equation}
respectively, and $\check{\bu}^{\caol,n+\hf}_i$
on the control volume $V_i$  is defined by
\[
\check{\bu}^{\caol,n+\hf}_i=-\frac{1}{\nu^\caol}\left[e^{\xi^{n+\hf}_i}\td\nabla_h Q^{\caol,n+\hf}_i+\td\nabla_h\left(z^\caol\psi^{n+\hf}+e^{\xi^{n+\hf}} \right)_i \right].
\]
Notice that $R^{n+\hf}$ is positive by the inequality $e^{-x}+\frac{e^{x}-e^{y}}{2e^{2x}}+\frac{(e^{x}-e^{y})^2}{3e^{3x}}>0$ for any $x, y \in \R$. 
 The initial and boundary conditions are discretized the same as in \reff{Bcs}. We refer to this scheme as ``{\bf Scheme II}" in what follows. 

Similar to the first-order scheme, we have
\begin{theorem}
The Scheme II \reff{DisPNPF2} enjoys the following properties:
\begin{compactenum}
\item[\rm (1)]Mass conservation: 
\[
\ciptwo{e^{\eta^{\caol,n+1}}}{1}=\ciptwo{e^{\eta^{\caol,n}}}{1}~~\mbox{for }~\caol=1,2,\cdots,M.
\]

\item[\rm (2)]Positivity preserving: Assume there exists a solution ($\eta^{\caol, n+1}$, $\psi^{n+1}$, $\xi^{n+1}$) to the Scheme II \reff{DisPNPF2}.  If $c^{\caol,n}_i>0$ and $T^n_i>0$, then $c^{\caol,n+1}_i>0$ and $T^{n+1}_i>0$ for $i=1,2,\cdots,N$.

\item[\rm (3)]
Discrete entropy increasing:\\
\end{compactenum}
\begin{equation}\label{ED2:eq0}
\quad\frac{S^{n+1}_h-S^n_h}{\dt}\geq \ep\sum_{\caol=1}^M \ciptwo{\nu^\caol e^{\eta^{\caol,n+\hf}} |\check{\bu}^{\caol,n+\hf}|^2}{R^{n+\hf}}
+ k [e^{\xi^{n+\hf}}\nabla_h \log R^{n+\hf},\nabla_h R^{n+\hf}]\geq 0.
\end{equation}
\end{theorem}

\noindent \proofc
(1) The proof of mass conservation is essentially the same as that of Theorem \ref{th:conservation}.

(2) The positivity follows from $c^{\caol, n+1}=e^{\eta^{\caol, n+1}}$ for $\caol=1,2,\cdots,M$ and $T^{n+1}=e^{\xi^{n+1}}$ with the assumption that there exists a solution ($\eta^{\caol, n+1}$, $\psi^{n+1}$, $\xi^{n+1}$)  to the Scheme II.

(3) We now proceed to prove the original entropy increase. Taking a discrete inner product of the fourth equation in \reff{DisPNPF2} with $\dt R^{n+\hf}$ yields
\begin{equation}\label{ED2:eq1}
\begin{aligned}
C_T\ciptwo{e^{\xi^{n+1}}-e^{\xi^{n}}}{R^{n+\hf}}=&-\dt k \ciptwo{\nabla_h\cdot(e^{\xi^{n+\hf}}\nabla_h \log R^{n+\hf})}{ R^{n+\hf}}\\
&+\dt \sum_{\caol=1}^M \left[\ep\ciptwo{\nabla_h\cdot (e^{\eta^{\caol,n+\hf}}\bu^{\caol,n+\hf} Q^{\caol,n+\hf})}{1}\right.\\
&\left.+\ciptwo{(1+Q^{\caol,n+\hf}) (e^{\eta^{\caol,n+1}}-e^{\eta^{\caol,n}})}{1}\right]\\
&+\dt\ep\sum_{\caol=1}^M \ciptwo{\nu^\caol e^{\eta^{\caol,n+\hf}} |\check{\bu}^{\caol,n+\hf}|^2}{R^{n+\hf}},
\end{aligned}
\end{equation}
where the equation for $P^{n+\hf}$ in \reff{PQR} has been used. For the first term on the right hand side, one can obtain
\begin{equation}\label{ED2:eq2}
\begin{aligned}
-\dt k\ciptwo{\nabla_h\cdot(e^{\xi^{n+\hf}}\nabla_h \log R^{n+\hf})}{ R^{n+\hf}}=&-\dt k\sum_{i=1}^N R^{n+\hf}_i\sum_{\sigma\in\partial V_i} \tau_\sigma e^{\xi^{n+1}_\sigma}D(\log R^{n+\hf})_{i,\sigma}\\
=&\dt k\sum_{\sigma\in\calE_{int}} \tau_\sigma e^{\xi^{n+1}_\sigma}D(\log R^{n+\hf})_{\sigma}D(R^{n+\hf})_{\sigma}\\
\geq& 0,
\end{aligned}
\end{equation}
where the thermally insulated boundary condition~\reff{Bcs} for $T$ has been used in the second equality, and monotonicity of the function $\log (x)$ has been used in the last inequality.
The following inequalities are available for the discrete inner products:
\begin{equation}\label{ED2:eq3}
\begin{aligned}
&\ciptwo{\nabla_h\cdot (e^{\eta^{\caol,n+\hf}}\bu^{\caol,n+\hf} Q^{\caol,n+\hf})}{1}=0,\\
&\sum_{\caol=1}^M\ciptwo{\nu^\caol e^{\eta^{\caol,n+\hf}} |\check{\bu}^{\caol,n+\hf}|^2}{R^{n+\hf}}\geq 0.
\end{aligned}
\end{equation}

The following Taylor expansion is valid: for $H(s)\in C^4(\R)$ and $x,~y\in\R$,
\[
\begin{aligned}
H(x)=&H(y)+H^{(1)}(y)(x-y)+\frac{1}{2} H^{(2)}(y)(x-y)^2+\frac{1}{6}H^{(3)}(y)(x-y)^3+\frac{1}{24}H^{(4)}(\theta)(x-y)^4,
\end{aligned}
\]
where $\theta$ is between $x$ and $y$, and $H^{(p)}(y)=\frac{\partial^{p}H}{\partial y^p}$ for $p=1,2,3,4$. 
Taking $H(s)=s\log s$, $x=e^{\eta^{\caol,n}_i}$, and $y=e^{\eta^{\caol,n+1}_i}$, one has $H^{(4)}=2s^{-3}>0$ for $s>0$. Therefore,
\[
e^{\eta^{\caol,n+1}_i}\eta^{\caol,n+1}_i-e^{\eta^{\caol,n}_i}\eta^{\caol,n}_i\leq (Q^{\caol,n+\hf}_i+1) (e^{\eta^{\caol,n+1}_i}-e^{\eta^{\caol,n}_i}).
\]
Taking $H(s)=\log s$, $x=e^{\xi^{n}_i}$, and $y=e^{\xi^{n+1}_i}$, one has $H^{(4)}=-6s^{-4}<0$ for $s>0$. Therefore, 
\[
\xi^{n+1}_i-\xi^n_i\geq R^{n+\hf}_i(e^{\xi^{n+1}_i}-e^{\xi^{n}_i}).
\]
Consequently, we have
\begin{equation}\label{ED2:eq4}
\ciptwo{e^{\xi^{n+1}}-e^{\xi^{n}}}{ R^{n+\hf} }\leq \ciptwo{\xi^{n+1}-\xi^n}{1},
\end{equation}
and
\begin{equation}\label{ED2:eq5}
\ciptwo{e^{\eta^{\caol,n+1}}}{ \eta^{\caol,n+1}}-\ciptwo{e^{\eta^{\caol,n}}}{ \eta^{\caol,n}}\leq \ciptwo{Q^{\caol,n+\hf}+1}{e^{\eta^{\caol,n+1}}-e^{\eta^{\caol,n}} }.
\end{equation}
Finally, a combination of \reff{ED2:eq1}--\reff{ED2:eq5} yields~\reff{ED2:eq0}. This completes the proof for Part (3).
\qed
\begin{remark}
The right hand side of \reff{ED1:eq0} and \reff{ED2:eq0} are first and second-order temporal discretization of the entropy production rate~\reff{dS/dt}, respectively.
\end{remark}

\begin{remark}
The existence of a solution ($\eta^{\caol, n+1}$, $\psi^{n+1}$, $\xi^{n+1}$) to the Scheme II is currently unavailable and will be studied in our future work.
\end{remark}

\section{Numerical Results}

\subsection{Accuracy Test}\label{s:test1}
We now test the accuracy of the Scheme I \reff{DisPNPF} and Scheme II \reff{DisPNPF2} for the PNPF system with two ionic species on a 2D computational domain $\Omega=[0,1]^2$. The nonlinear numerical schemes are solved iteratively with the Newton's iterations 
 at each time step. Consider the following dimensionless PNPF system:
\begin{equation}\label{eq:ex}
\left\{
\begin{aligned}
&\partial_t c^\caol+\nabla\cdot (c^\caol \bu^\caol)=f_\caol,\\
& c^\caol \bu^\caol=-c^\caol T\nabla \log c^\caol - c^\caol \nabla(\psi+T),~ \caol=1, 2,\\
&-\Delta\psi=c^1-c^2+\rho^f,\\
&\partial_t T=\Delta T+ \sum_{\caol=1}^2 \left\{T\left[\nabla\cdot(c^\caol\bu^\caol\log c^\caol)+(1+\log c^\caol)\partial_t c^\caol \right]+ c^\caol |\bu^\caol |^2\right\}+f_3,
\end{aligned}
\right.
\end{equation}
The source terms $f_1$, $f_2$, $\rho^f$, and $f_3$ are determined by the following exact solution
\begin{equation}\label{ExS}
\left\{
\begin{aligned}
&c^1=0.1e^{-t}\cos(\pi x)\cos(\pi y)+0.2,\\
&c^2=0.1e^{-t}\cos(\pi x)\cos(\pi y)+0.2,\\
&T=0.1e^{-t}\cos(\pi x)\cos(\pi y)+0.2,\\
&\psi=0.1e^{-t}\cos(\pi x)\cos(\pi y).
\end{aligned}
\right.
\end{equation}
The initial conditions are obtained by evaluating the exact solution at $t=0$. We consider zero-flux boundary conditions for concentrations, thermally insulated boundary conditions for temperature,
and the following boundary conditions for electric potential:
\[
\left\{
\begin{aligned}
 &\psi(t,0,y)=\frac{1}{10}e^{-t}\cos(\pi y),~\psi(t,1,y)=-\frac{1}{10}e^{-t}\cos(\pi y), \quad y\in[0,1],\\
&\frac{\partial\psi}{\partial y}(t,x,0)=\frac{\partial\psi}{\partial y}(t,x,1)=0,~~x\in [0,1].\\
\end{aligned}
\right.
\]

\begin{figure}[htbp]
\centering
\subfigure[Scheme I]{\includegraphics[scale=.36]{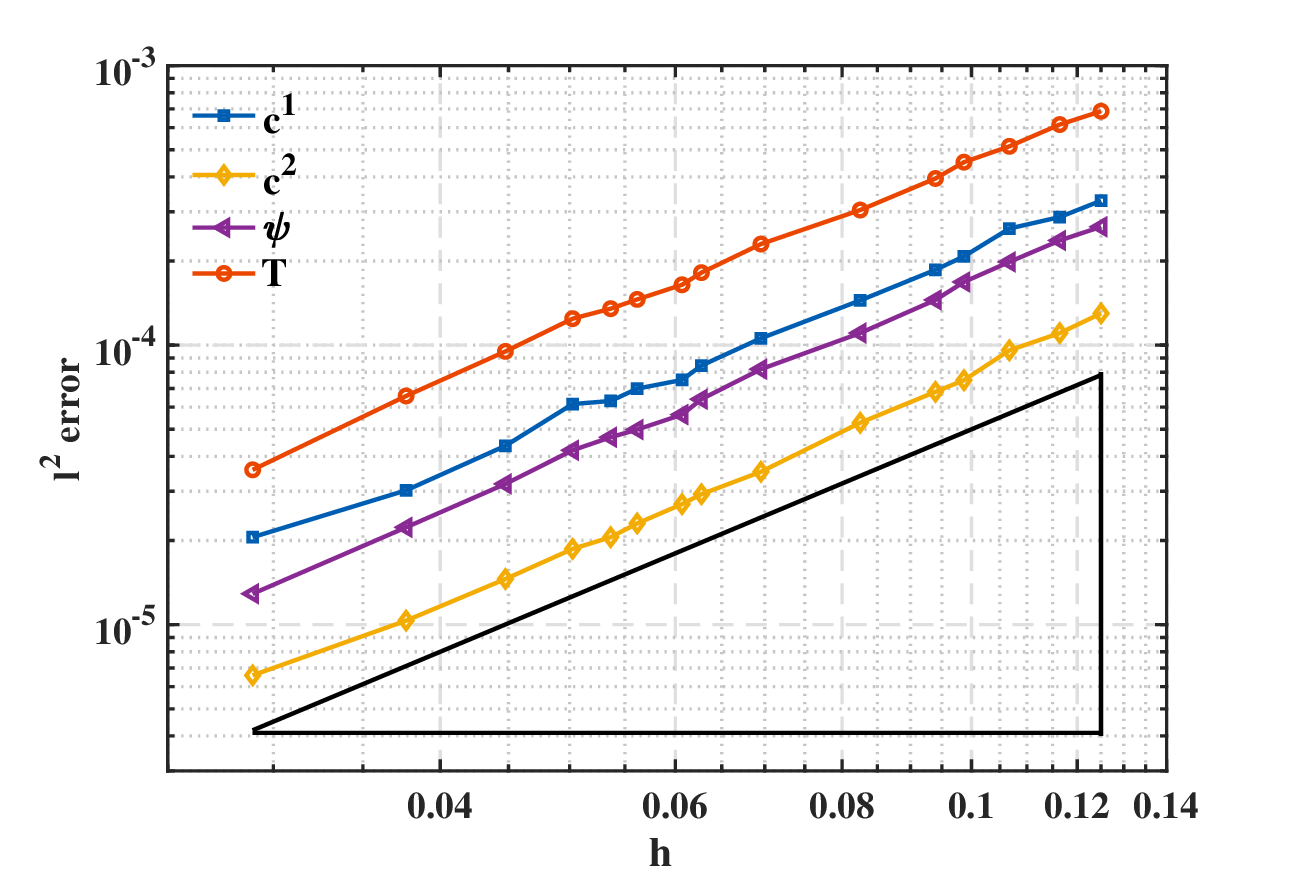}}
\subfigure[Scheme II]{\includegraphics[scale=.36]{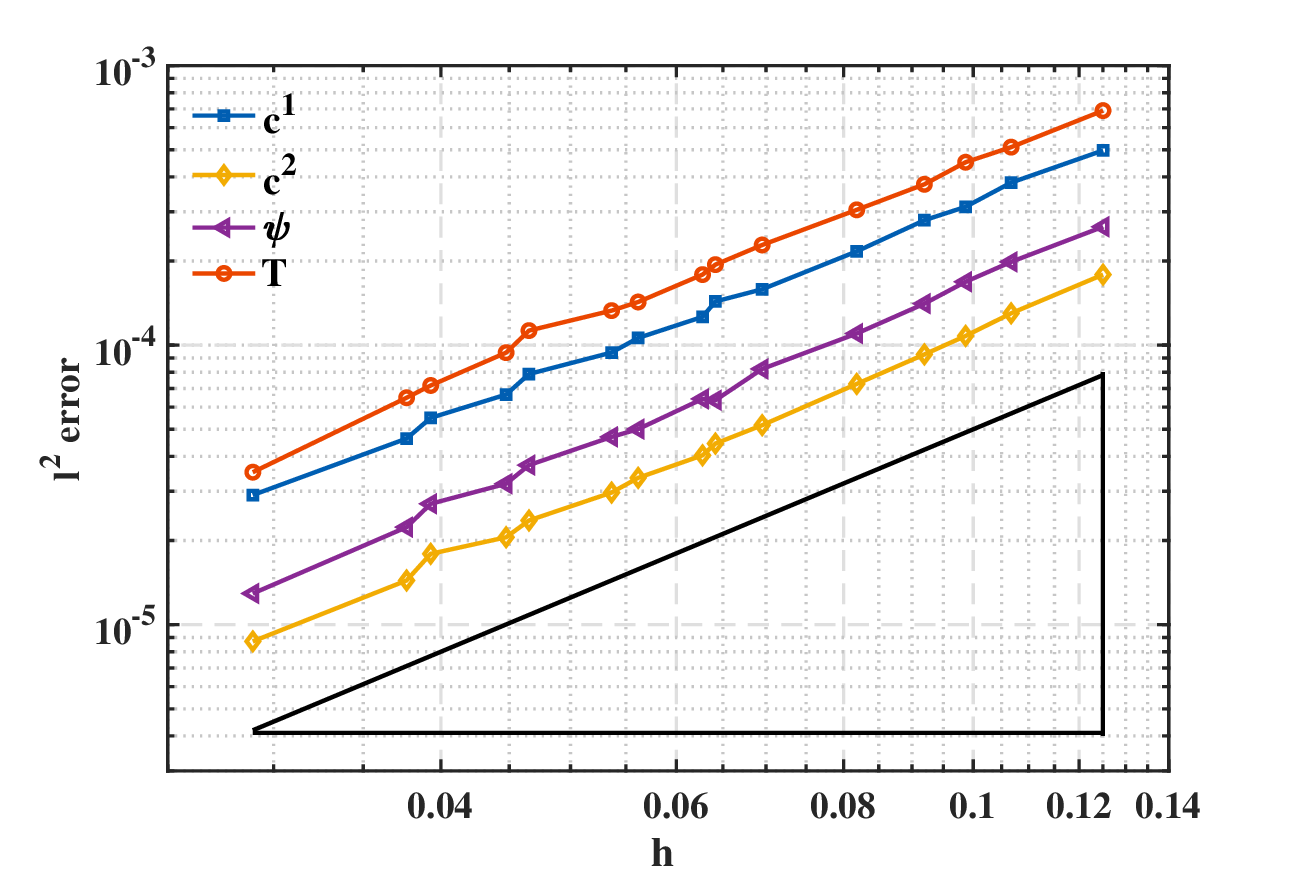}}
\caption{Numerical error of $c^1$, $c^2$, $\psi$, and $T$ at time $T=0.1$ obtained by (a) Scheme I with a mesh ratio $\Delta t=h^2$ and (b) Scheme II with a mesh ratio $\Delta t=h/10$.}
\label{f:error1}
\end{figure}
We first test the numerical accuracy of Scheme I utilizing various spatial step sizes $h$ with a fixed mesh ratio $\Delta t=h^2$. Figure~\ref{f:error1} (a) records discrete $l^2$-error of ionic concentration, electrostatic potential, and temperature at time $T=0.1$.  One can observe that the error decreases as the mesh refines. Comparison with the reference slope implies that the convergence rate for both ion concentrations, electrostatic potential and temperature approaches $O(h^2)$ as $h$ decreases. This indicates that the Scheme I, as expected, is first-order and second-order accurate in time and spatial discretization, respectively. Note that the mesh ratio here is chosen for the purpose of accuracy test, not for stability or positivity.

Next, we test the numerical accuracy of Scheme II with a mesh ratio  $\Delta t=h/10$. As displayed in Figure~\ref{f:error1} (b), the numerical error decreases with a convergence order around 2 as well, indicating that Scheme II \reff{DisPNPF2} is second-order in both time and spatial discretization.

\subsection{Property Tests}
\begin{figure}[H]
\centering
\includegraphics[scale=.4]{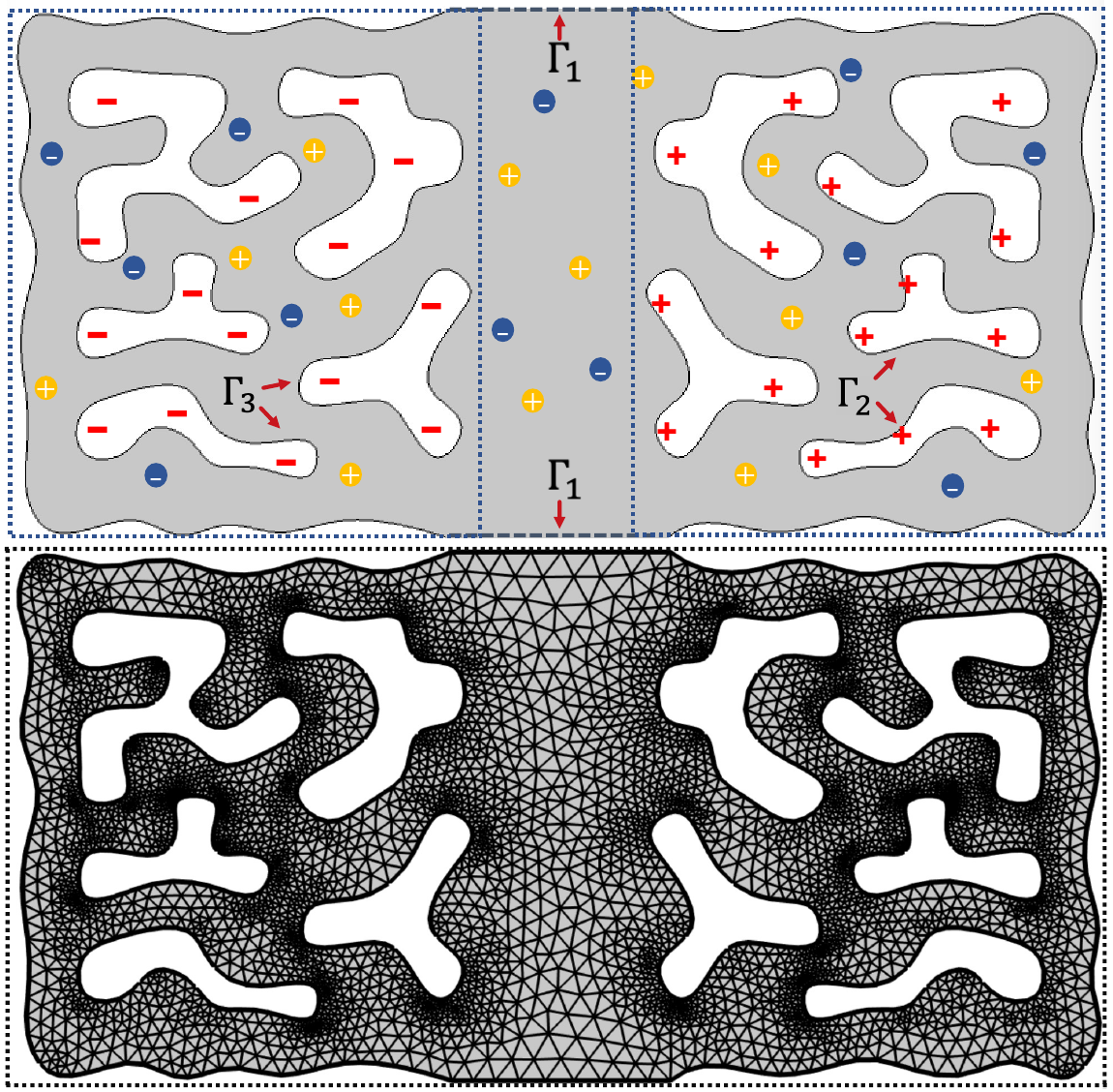}
\caption{A schematic illustration of the computational domain $\Omega$ filled with electrolytes shown in gray (upper) and the mesh (lower) contained within a rectangle $[-10,10]\times[0,10] ~{\rm nm}^2$. Biased voltage differences are applied across two thermally-insulated, blocking electrodes with complex interfaces denoted by $\Gamma_2$ and $\Gamma_3$. The boundary in the middle is labeled by $\Gamma_1$. 
}
\label{f:domain}
\end{figure}
In this section, we demonstrate the performance of the proposed numerical schemes in preserving properties of physical significance at discrete level. We numerically explore the ionic and heat dynamics in a supercapacitor with electrodes of complex geometry, as shown in Figure~\ref{f:domain}.  The boundary of the computational domain, $\partial\Omega$, consists of three disjoint parts: 
\[
\begin{aligned}
&\Gamma_1=\left\{(x,y): y=0~{\rm nm}\mbox{ or } 10~{\rm nm},~ -2 ~{\rm nm}\leq x\leq 2 ~{\rm nm} \right\},\\
&\Gamma_2=\left\{(x,y): (x,y)\in\partial\Omega \backslash \Gamma_1,~ x>2 ~{\rm nm} \right\},\\
&\Gamma_3=\left\{(x,y): (x,y)\in\partial\Omega \backslash \Gamma_1,~ x<-2 ~{\rm nm}\right\}.
\end{aligned}
\]   
Note that $\Gamma_2\cup\Gamma_3=\partial\Omega\backslash\Gamma_1$.
To study ion and heat transport through electrodes of complex geometry under an applied voltage, we consider binary monovalent electrolytes and set initial conditions as 
\[
c^1(0,x,y) = 0.2~{\rm M},~ c^2(0,x,y)=0.2~{\rm M},~T(0,x,y)=300~{\rm K} \mbox{~for } (x,y)\in \Omega.
\]
Also, we prescribe zero-flux boundary conditions for ionic concentrations, thermally insulated boundary conditions for temperature, and the following boundary conditions for electric potential:
\begin{equation}\label{InitBouCons2}
\left\{
\begin{aligned}
&\psi(t,x,y)=\psi_* ~~{\rm k_B T_0/e}~~ &&\text{on } \Gamma_2,\\
&\psi(t,x,y)=0 ~~{\rm k_B T_0/e}~~ &&\text{on } \Gamma_3,\\
&\nabla\psi\cdot \n=0 ~~{\rm k_B T_0/(e\mu m)}~~ &&\text{on } \Gamma_1,\\
\end{aligned}
\right.
\end{equation}
which describes a horizontally applied voltage, and zero surface charge on upper and lower boundaries. 
Unless otherwise stated, the following numerical simulations take $\psi_* =2 ~{\rm k_B T_0/e}$, $\ve_0=8.85\times 10^{-21}~{\rm C/(Vnm)}$, $\ve_r=80$, $C=38.8 $ M, $k=1.20\times 10^{-4} {\rm J/(K ms)}$, $\nu^1=\nu^2=4.14\times 10^{-10} {\rm (Js)/m^2}$, $z^1=1$, $z^2=-1$, and $\rho^f=0~{\rm M}$. 

\begin{figure}[H]
\centering
\includegraphics[scale=.43]{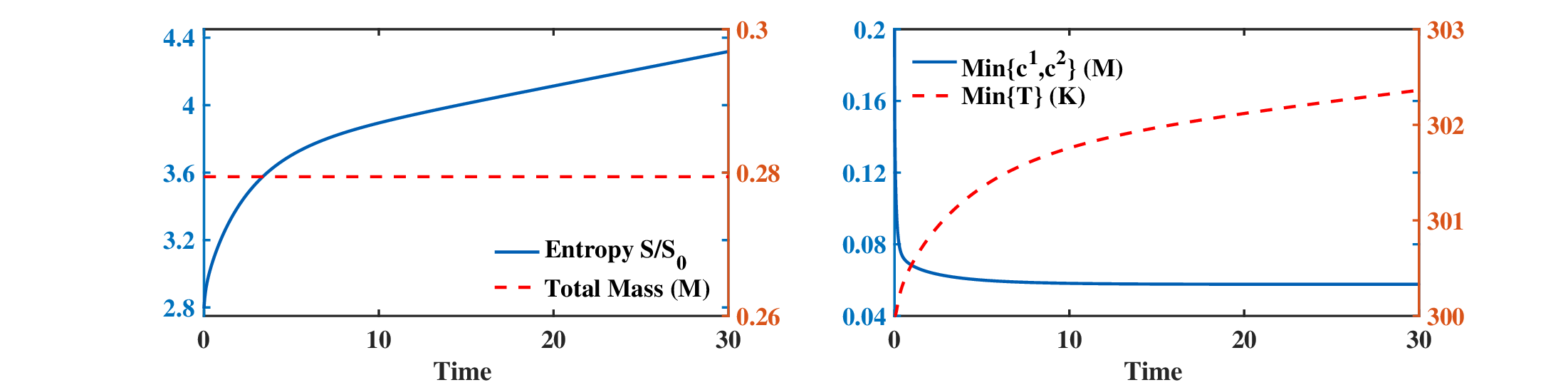}
\caption{Evolution of the discrete entropy $S$ (rescaled by $S_0:=k_Bc_0L^2$) and total mass of cations (Left), as well as minimum concentration and temperature over the computational mesh (Right).}
\label{f:MassEntropy}
\end{figure}
With zero-flux mass, thermally insulated,  and time-independent voltage boundary conditions, the system possesses properties of mass conservation and positive entropy production. Figure~\ref{f:MassEntropy} presents the evolution of discrete entropy $S_h^n$ (cf.~\reff{DisEntropy}) and total mass of cations, as well as the minimum concentration ${\rm Min}_i \left\{c^1_{i}, c^2_{i}\right\}$ and minimum temperature ${\rm Min}_i \left\{T_{i}\right\}$.
From the left panel of Figure~\ref{f:MassEntropy}, one observes that the entropy \reff{DisEntropy} increases monotonically and total mass remains constant as time evolves. The right panel of Figure~\ref{f:MassEntropy} demonstrates that the minimum of cations and $T$ on the computational mesh maintain positive, indicating that the developed numerical schemes preserve positivity at discrete level. Such numerical results further confirm our analysis on property preservation.
\begin{figure}[htbp]
\centering
\includegraphics[scale=.51]{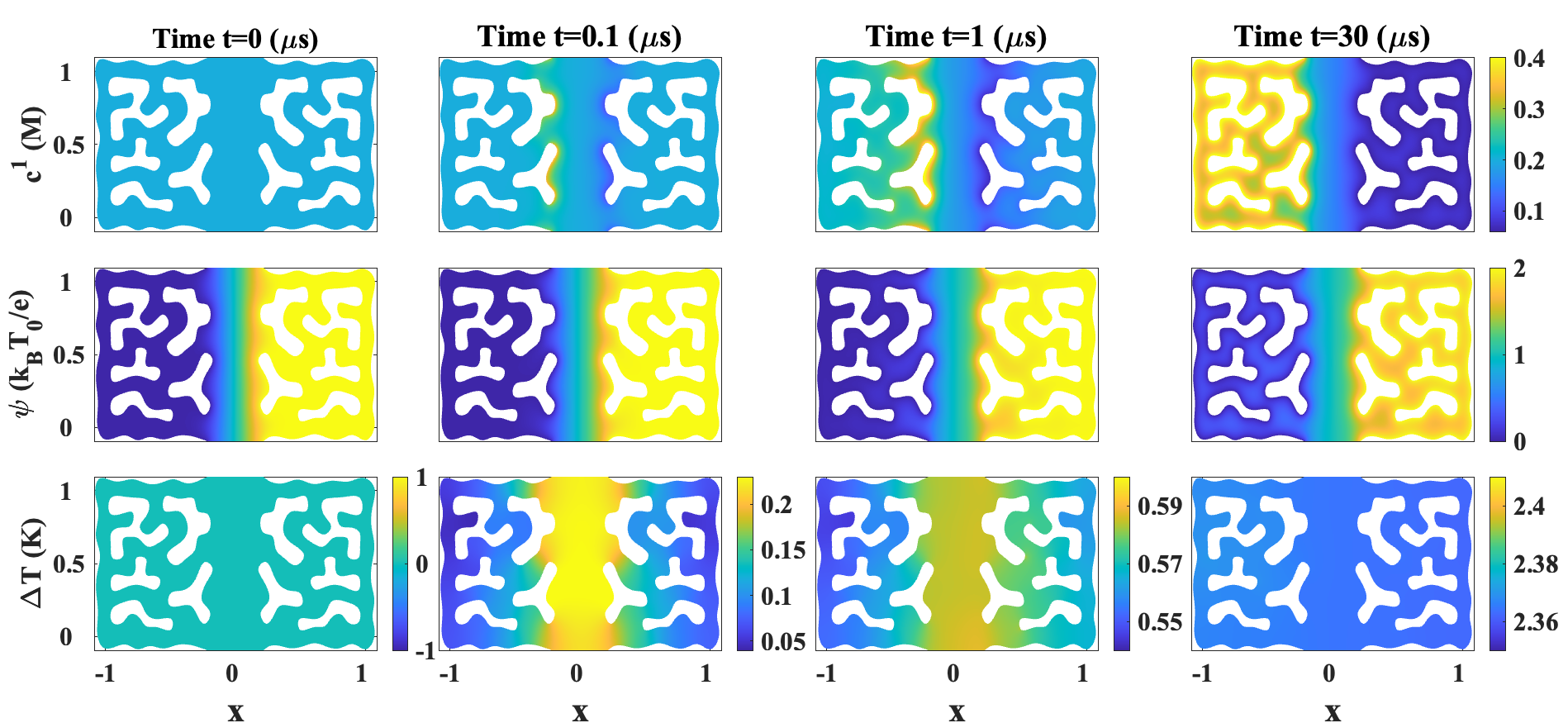}
\caption{Evolution of the cation concentration $c^1$, the electric potential $\psi$, and temperature rise $\Delta T:=T-T_0$ at time $t=0 $ $t=0.1$, $t=1$ and $t=30 ({\rm \mu s})$.}
\label{f:Evolution}
\end{figure}

\subsection{Numerical Investigation}
\subsubsection{Charging Process}
We next apply our PNPF model, along with the proposed numerical methods, to study non-isothermal electrokinetics in supercapacitors with  electrodes of complex geometry (cf. Figure~\ref{f:domain}), when it is charged towards a steady state. A fixed voltage difference, $\psi_* = 2 ({\rm k_B T_0/e})$, is applied across electrodes. Figure~\ref{f:Evolution} displays the evolution of cation concentration $c^1$, electric potential $\psi$, and temperature $T$ at time $t=0$, $t=0.1$, $t=1$, and $t=30\,{\rm \mu s}$.  As the charging proceeds, the cations $c^1$, as counterions to the left electrode, gradually permeate into the left electrode and accumulate next to the irregular boundary, forming electric double layers. Meanwhile, the cations also get depleted away from the right electrode due to electrostatic interactions. On the other hand, the electric potential, $\psi$, gets screened gently by the permeated counterions in both the cathode and anode.  

We next discuss the temperature distribution and evolution during the charging process. The lower plot of Figure~\ref{f:Evolution} displays the dynamics of temperature rise $\Delta T(x, t):=T-T_0$. At time $t = 0.1$, the temperature first starts to rise most obviously in the central bulk region due to the large current term $\sum_{\caol=1}^M \nu^\caol c^\caol |\bu^\caol |^2$ across the cathode and anode, which is always positive. Analogous to the Joule heating, such large convection is responsible for the main heat generation in initial fast timescale of the charging process.
At time $t = 1$, the high temperature in the central bulk region starts to diffuse into the branching region of electrodes, accompanying the permeation of ions into electrodes in a slower timescale. The spatial temperature heterogeneity gradually gets smoothed out in the later stage of the charging process. As approaching the steady state, the temperature in the supercapacitor becomes homogeneous with a rise of temperature around $2.37$~K.
\begin{figure}[htbp]
\centering
\subfigure{\includegraphics[scale=.38]{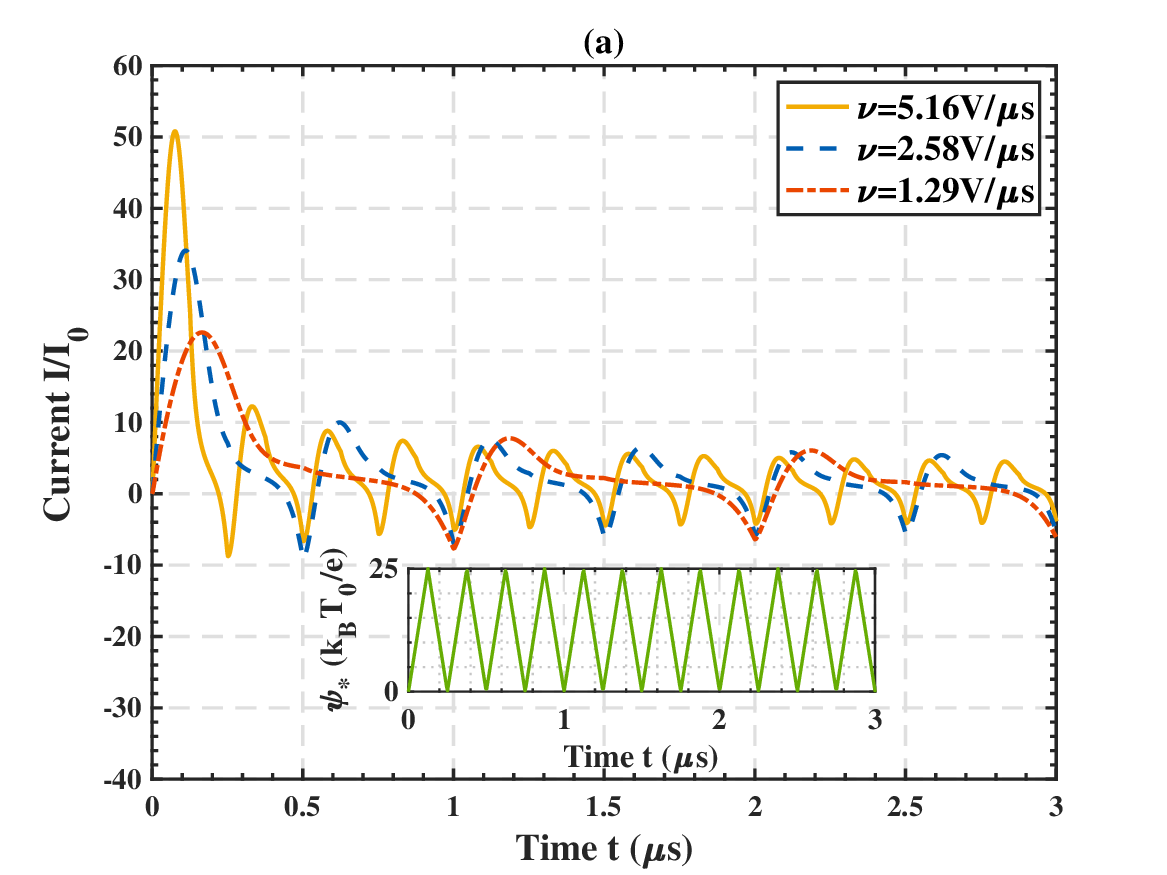}}\hspace{0.01\linewidth}
\subfigure{\includegraphics[scale=.38]{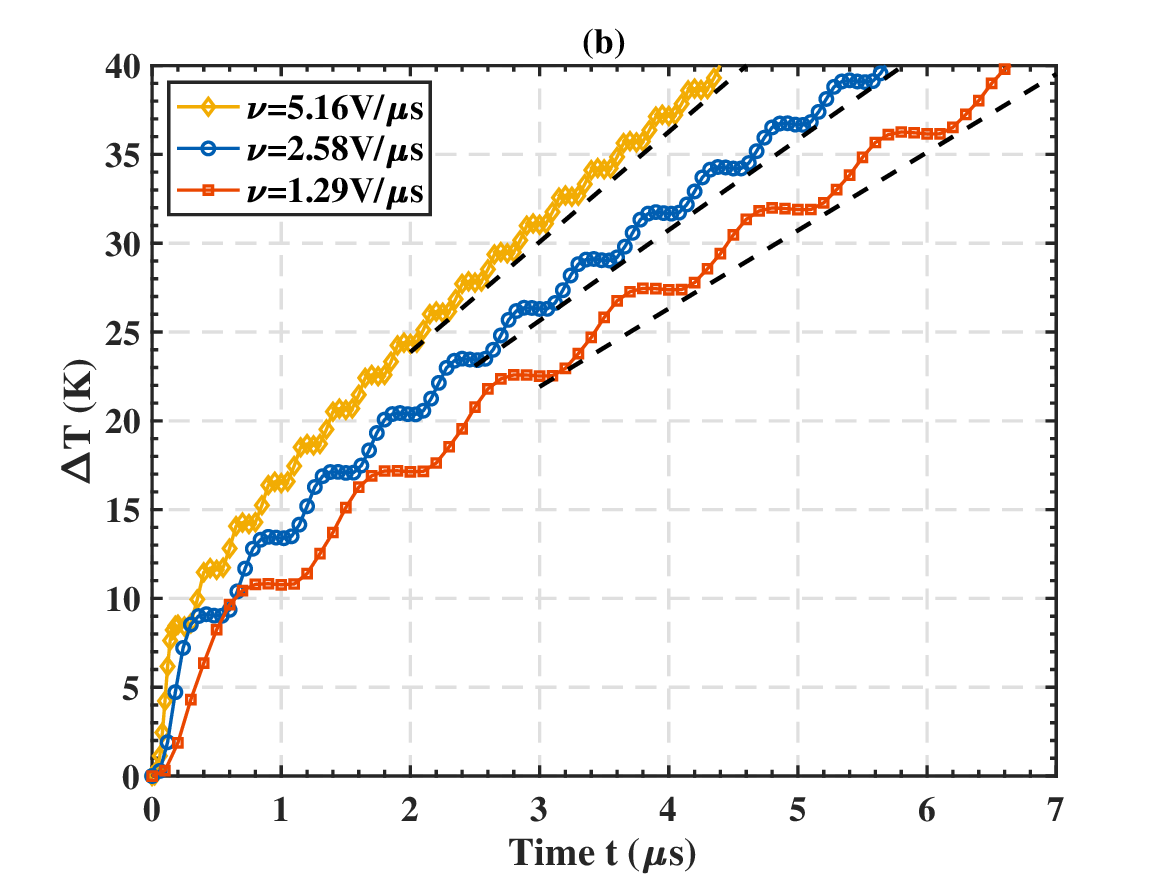}}\hspace{0.01\linewidth}
\subfigure{\includegraphics[scale=.38]{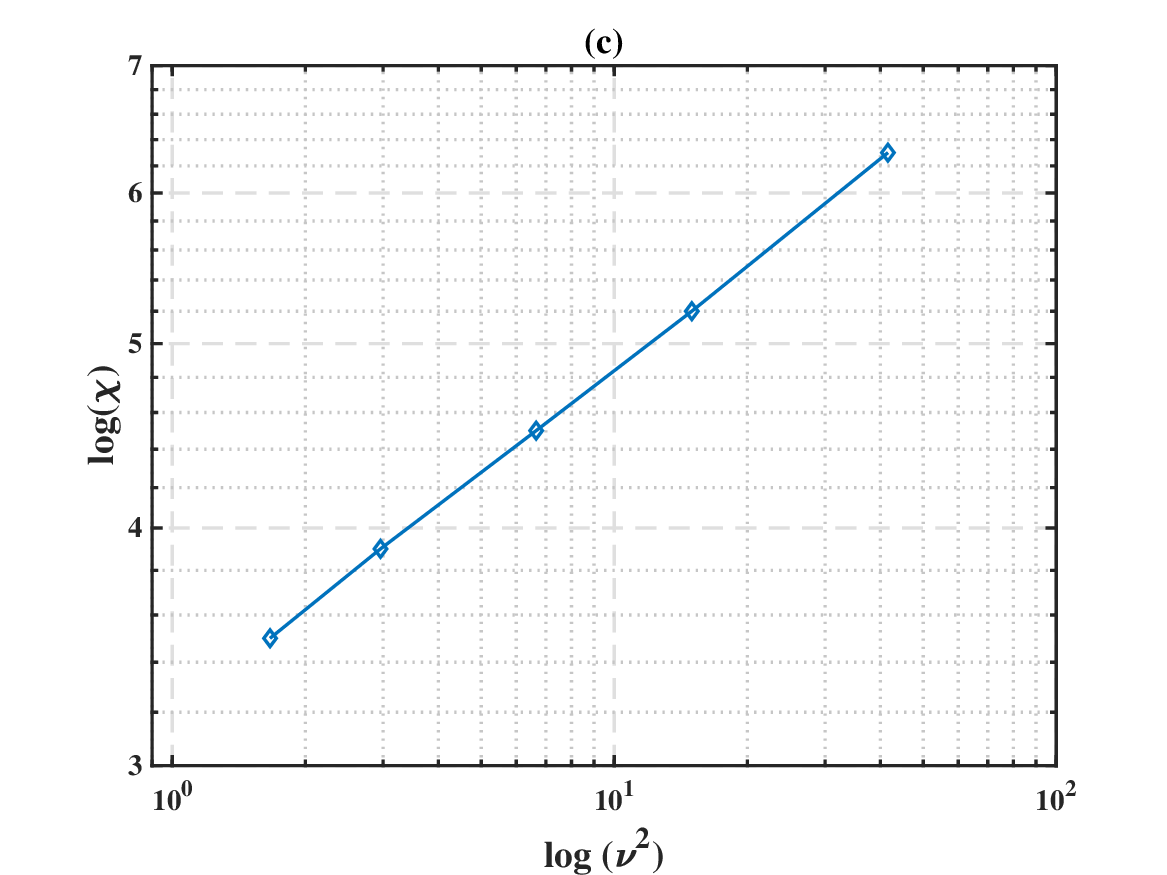}}
\subfigure{\includegraphics[scale=.38]{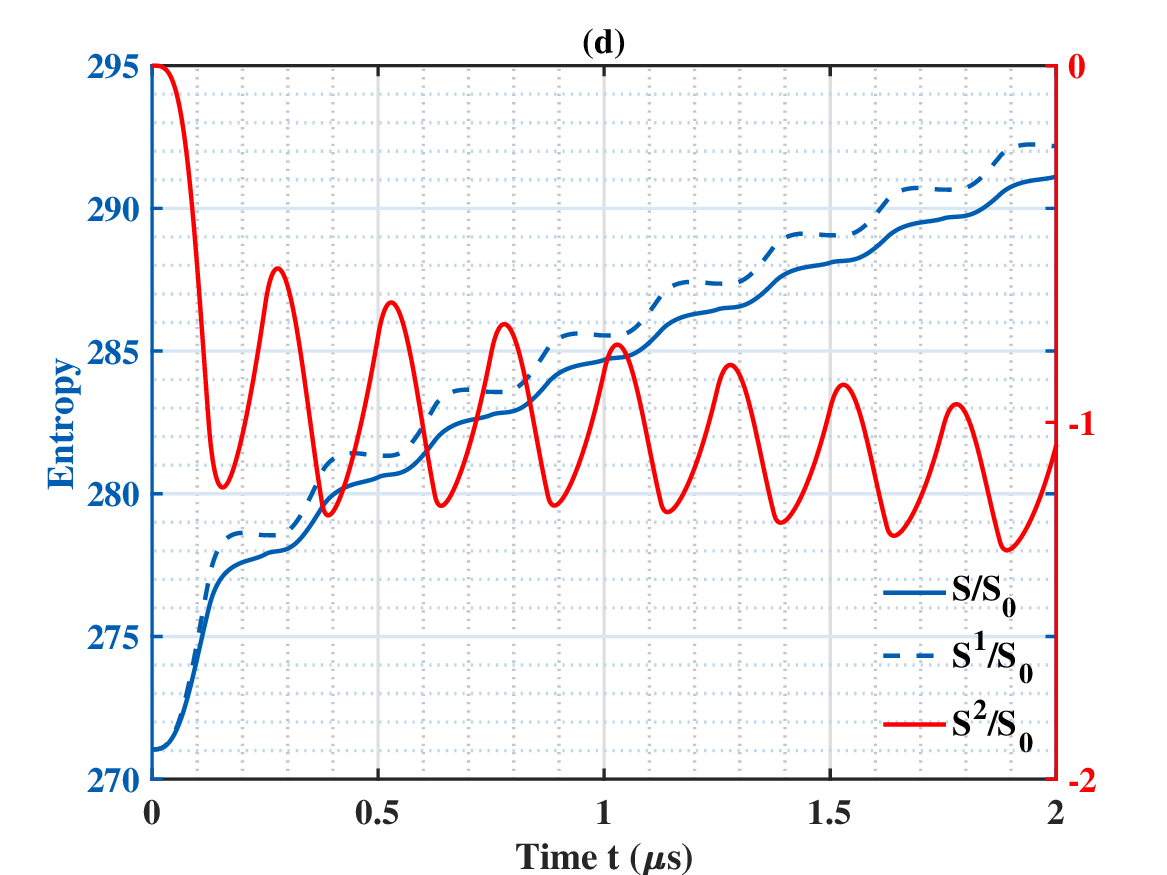}}
\caption{(a) The evolution of current density $I$ (rescaled by $I_0:=k_BT_0 c_0/ \nu_0 L$) in the bulk, as well as a linearly varying scheme of an applied voltage in the inset; (b) Temperature change $\Delta T$ using different scan rate $\nu$; (c) A log-log plot of temperature rising slope $\chi$ against $\nu^2$, showing a quadratic scaling law; (d) Entropy $S_h$ (rescaled by $S_0$) and its two components. }
\label{f:T-Time}
\end{figure}
\subsubsection{Cyclic Voltammetry}
Cyclic voltammetry (CV) is a popular and powerful electrochemical technique commonly used to evaluate the performance of supercapacitors~\cite{WangPilon_EA12,WangPilon_JPCC13,GirardWangPilon_JPCC15,MeiPilon_EA17}. We apply the proposed numerical methods to probe non-isothermal electrokinetics in supercapacitors with complex electrode geometry using a linearly varying applied voltage $\psi_*$ with the following scheme
\begin{equation}\label{cvPotential}
\psi_*(t)=\left\{
\begin{aligned}
&\nu t,~~&&2(n-1)t_0\leq t\leq (2n-1)t_0,\\
&25-\nu \left[t-(2n-1)t_0 \right],~~&&(2n-1)t_0\leq t\leq 2nt_0,\\
\end{aligned}
\right.
\end{equation}
where $\nu$ is the scan rate in $V/\mu s$ and $n (=1,2,3,\cdots)$ is the number of cycles; cf. the inset plot of Figure \ref{f:T-Time} (a). 

Under a periodically oscillatory applied voltage, one can see from Figure~\ref{f:T-Time} (a) that the current rises drastically in the first charging period and dives into negative values in the first discharging period. Later on, the current rises and decreases in a stable, periodic manner in charging/discharging cycles.  Similar response of temperature rise can be seen in Figure~\ref{f:T-Time} (b), in which the temperature is heated up monotonically in a charging stage and slightly cooled down in a discharging stage. Overall, one can see that the temperature increases in an oscillatory way, indicating that our approach can robustly capture the reversible and irreversible heating in CV tests. In addition, one can find that temperature rising speed depends directly on the scanning rate $\nu$. To further quantify the temperature rise, we denote by $\chi$ the slope of overall temperature rise. Figure \ref{f:T-Time} (c) presents a log-log plot of the relation between $\chi$ and $\nu^2$, showing a perfect quadratic scaling law $\chi \propto \nu^2$.  Such a conclusion is consistent with previous understanding on Joule heating effect~\cite{Pilon_JPS15, Janssen_PRL17, JiLiuLiuZhou_JPS2022}.

Attention should be paid to the variation of the entropy and its two components: $S_1$ and $S_2$ (cf.~\reff{S1+S2}), in the Figure~\ref{f:T-Time} (d). As expected, the total entropy keeps increasing monotonically and oscillatory throughout the charging/discharging processes, since zero-flux and thermally insulated boundary conditions have been considered for the system. The entropy related to temperature, $S_1$, increases in an oscillatory manner due to the temperature variation. Notably, the entropy related to the randomness of ionic distribution, $S_2$, decreases in the charging stage, which can be attributed to the adsorption of counterions into the EDL during charging, resulting in a more ordered double-layer structure.  On the contrary, $S_2$ increases in the discharging stage, due to the dissolving ions back into the bulk region, leading to a more disordered state. Such results align with previous studies that entropy related to the ionic distribution accounts for the reversible heating in the charging/discharging processes through the formation and dissolution of EDLs~\cite{Pilon_JPS15, JiLiuLiuZhou_JPS2022}.

\begin{figure}[htbp]
	\centering
	\subfigure{\includegraphics[scale=.4]{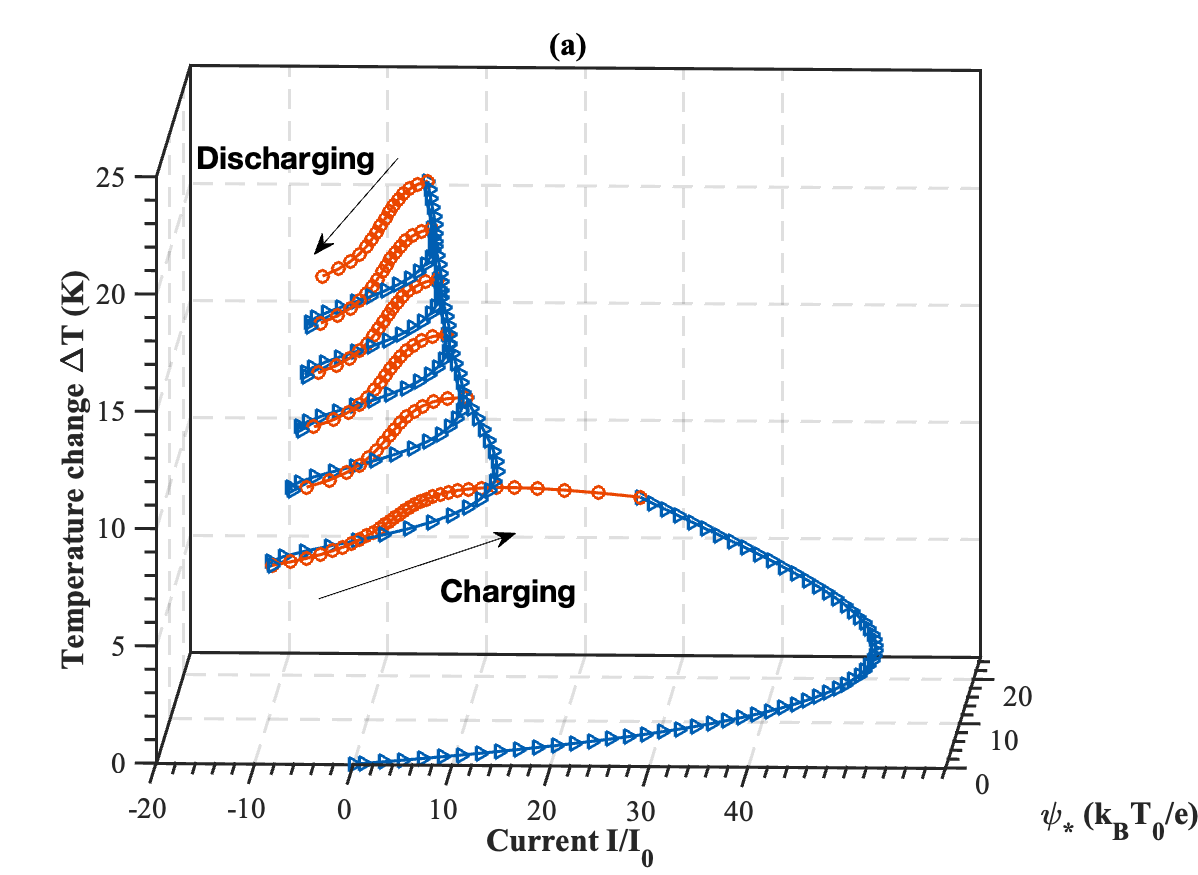}}
	\subfigure{\includegraphics[scale=.4]{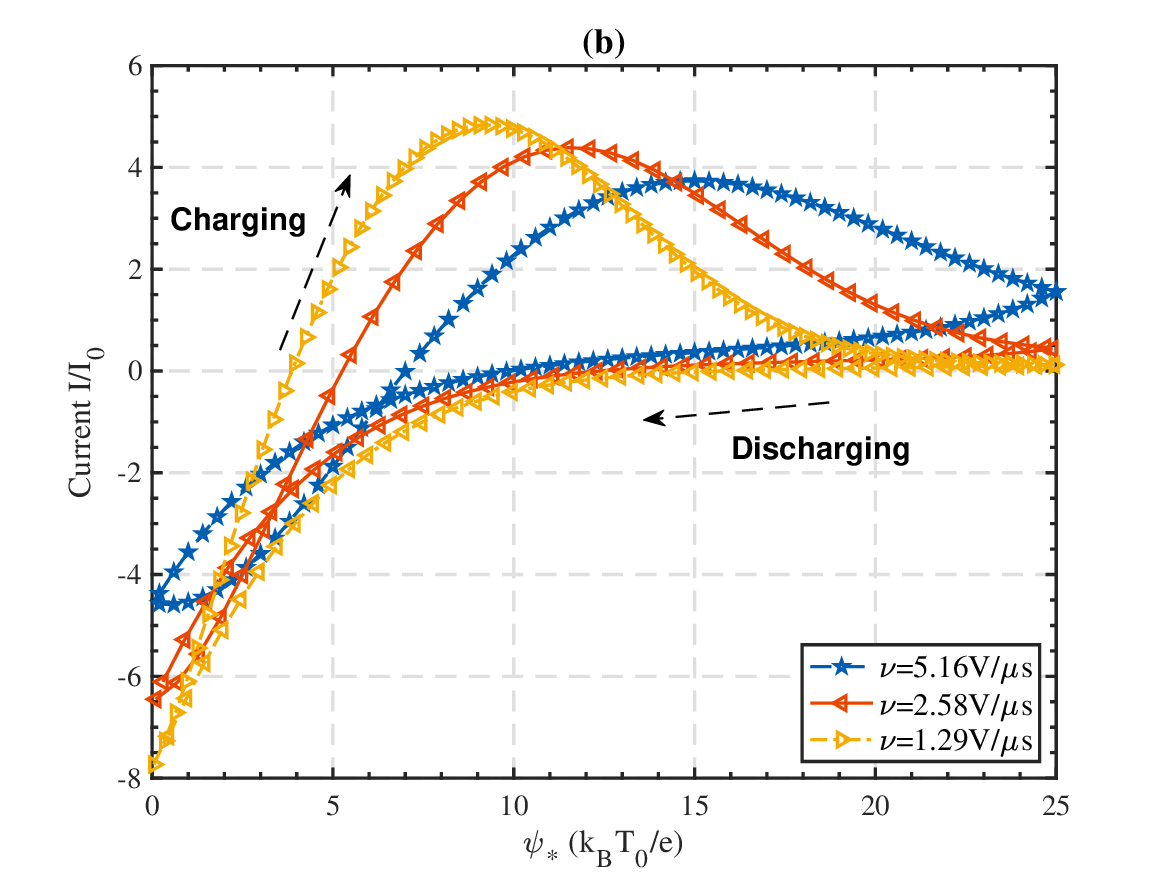}}
	\caption{(a) Evolution curve for the current in the bulk $I(0,t)/I_0$, the applied voltage $\psi_*(t)$, and the temperature change $\Delta T(t)$ in a CV test; (b) The curve of $I/I_0$ versus the applied voltage $\psi_*(t)$ with various scanning rates.}
	\label{f:IV}
\end{figure}
We apply our approach to understand the interplay between temperature and other key factors in CV measurements. Again, a linearly varying surface potential $\psi_*$~\reff{cvPotential} is imposed on electrodes of complex geometry. Figure \ref{f:IV}(a) presents a 3D evolution curve of the charge current in the bulk $I(0,t)/I_0$, the applied voltage $\psi_*(t)$, and the temperature change $\Delta T(t)$ in several charging/discharging cycles. In contrast to the traditional plot of current-versus-potential ($I$-$V$ curve), such a $I$-$V$-$T$ presentation includes the temperature as another dimension and unravels the dependence of the temperature on both the charge current and applied voltage. Overall, one can observe that the temperature rises in a spiral path, with rising temperature in charging and decreasing temperature  . It is of great interest to see that the area enclosed by the projected $I$-$V$ loop shrinks as the temperature rises, indicating that the thermal effect on electrokinetics, i.e., Soret effect, is faithfully captured by our model and numerical methods. To further probe the CV measurements, we perform simulations with different scanning rates. Figure~\ref{f:IV} (b) displays the $I$-$V$ diagram with scanning rates ranging from $\nu=1.29V/\mu s$ to $\nu=5.16V/\mu s$. It is worth noting that the charging/discharging processes follow different paths, showing a hysteresis diagram. With faster scanning rates, one finds that there is an intersection point in the hysteresis curve, which is consistent with recent CV measurements in electrodes with complex pore geometry~\cite{2024_Bisquert_PRXEnergy}. However, the intersection point is missing in CV simulations for planar electrodes for which the problem can be reduced to a $1$D case~\cite{JiLiuLiuZhou_JPS2022}, highlighting the remarkable difference between planar electrodes and electrodes with complex geometry. This also emphasizes the significance of simulating electrodes with complex geometry for realistic supercapacitors.

\section{Conclusions}
In this work, finite-volume schemes have been proposed for a thermodynamically consistent PNPF model for the prediction of non-isothermal electrokinetics in supercapacitors. Semi-implicit discretization has been employed to develop the first-order accurate scheme, in which the decoupling of temperature from ionic concentrations and electric potential makes the scheme efficient to solve.  Novel modified Crank-Nicolson discretization of $\log \frac{1}{T}$ and $\log c^{\caol}$ has been designed to develop a second-order accurate scheme. Rigorous numerical analysis on the first-order scheme has established the unique existence of positive ionic concentrations and temperature at discrete level. Moreover, both the first-order and second-order schemes have been proved to unconditionally respect the mass conservation and original entropy increase at discrete level.  Numerical experiments have been conducted to verify that our numerical methods have expected accuracy and are capable of preserving anticipated properties.  Simulation results have demonstrated that our PNPF model, along with the designed numerical schemes, can successfully predict temperature oscillation in the charging/discharging processes of supercapacitors with electrodes of complex geometry. Furthermore, temperature rising slope is found to scale quadratically against the scanning rate in CV tests. In addition to the monotonically rising total entropy, simulations have unravel that the ionic entropy contribution, which measures the disorder of ionic distribution, decreases in the formation of EDLs and increases in the dissolution of EDLs, highlighting that our approach is indeed able to capture underlying mechanisms responsible for reversible heat generation in the charging/discharging processes of supercapacitors. In summary, our work provides a promising tool in predicting non-isothermal electrokinetics in charging/discharging processes of supercapacitors with electrodes of complex geometry.

\vskip 5mm
\noindent{\bf Acknowledgements.}
The authors would like to thank Professor Chun Liu for helpful discussions. This work is supported in part by the National Natural Science Foundation of China 12101264, Natural Science Foundation of Jiangsu Province BK20210443, High level personnel project of Jiangsu Province 1142024031211190 (J. Ding), and National Natural Science Foundation of China 12171319 (X. Ji. and S. Zhou).

\bibliographystyle{plain}
\bibliography{PNP}

\end{document}